\documentclass[final]{siamltex}

\usepackage{graphicx}
\usepackage{gensymb}
\usepackage{amsmath}
\usepackage{makecell}
\usepackage{epsfig}
\usepackage{amsmath, amssymb}
\usepackage {graphicx}
\usepackage{float}
\usepackage{amsmath, amssymb}
\usepackage{mathptmx}      
\usepackage{multirow}
\usepackage{amssymb}
\usepackage{amsfonts}
\usepackage{mathrsfs}
\usepackage{amsmath, amssymb}
\usepackage{array}
\usepackage{subfigure}
\usepackage[graphicx]{realboxes}

\newtheorem{Lemma}{Lemma}[section]

\newtheorem{remark}{Remark}[section]
\newtheorem{algorithm}{Algorithm}

\title{A Semi-Randomized and Augmented Kaczmarz Method with Simple Random Sampling for Large-Scale Inconsistent Linear Systems}

\author{Shunchang Li\thanks{School of Mathematics, China University of Mining and Technology, Xuzhou 221116, Jiangsu, P.R. China. E-mail: {\tt lishunchang2021@163.com.} }
	\and Gang Wu\thanks{Corresponding author. School of Mathematics,
		China University of Mining and Technology, Xuzhou 221116, Jiangsu, P.R. China.
		E-mail: {\tt gangwu@cumt.edu.cn}. This work is supported by the National Natural Science Foundation
of China under Grant 12271518, the Fujian Natural Science Foundation under Grant 2023J01354, the Key
Research and Development Project of Xuzhou Natural Science Foundation under Grant KC22288, and the Open
Project of Key Laboratory of Data Science and Intelligence Education of the Ministry of Education under
Grant DSIE202203.}}
\begin{document}
\maketitle

\begin{abstract}
A greedy randomized augmented Kaczmarz (GRAK) method was proposed in [{\sc Z.-Z. Bai and W.-T. Wu}, SIAM J. Sci. Comput., 43 (2021), pp. A3892--A3911] for large and sparse inconsistent linear systems. However, one has to construct two new index sets via computing residual vector with respect to the augmented linear system in each iteration. Thus, the computational overhead of this method is large for extremely large-scale problems. Moreover, there is no reliable stopping criterion for this method. In this work, we are interested in solving large-scale sparse or dense inconsistent linear systems, and try  to enhance the numerical performance of the GRAK method. First, we propose an accelerated greedy randomized augmented Kaczmarz method. Theoretical analysis indicates that it converges faster than the GRAK method under very weak assumptions. Second, in order to further release the overhead, we propose a semi-randomized augmented Kaczmarz method with simple random sampling. Third, to the best of our knowledge, there are no practical stopping criteria for all the randomized Kaczmarz-type methods till now. To fill-in this gap, we introduce a practical stopping criterion for Kaczmarz-type methods, and show its rationality from a theoretical point of view. Numerical experiments are performed on both real-world and synthetic data sets, which demonstrate the efficiency of the proposed methods and the effectiveness of our stopping criterion.
\end{abstract}

\pagestyle{myheadings} \thispagestyle{plain} \markboth{S. LI AND G. WU}{\sc Semi-Randomized and Augmented Kaczmarz Method}

\date{ }%

\section{Introduction}
\setcounter{equation}{0}

Consider the large-scale linear system
\begin{equation}\label{1.1}
A{\bf x}={\bf b},
\end{equation}
where $A\in \mathbb{C}^{m\times n}$, ${\bf b}\in \mathbb{C}^{m}$, and ${\bf x}\in \mathbb{C} ^n$ is an unknown vector. In this paper, we make the assumption that there are no zero rows or columns in $A$. On one hand, if the linear system \eqref{1.1} is consistent, the interest is to find the least-norm solution ${\bf x}_{\star}=A^{\dagger}{\bf b}$.
The Kaczmarz method \cite{ref1} is an effective iterative method for consistent linear systems, which has been applied in many real-world applications, e.g., signal processing \cite{ref15, ref16, ref17}, image and sound reconstructions \cite{ref46, ref47, ref18, ref19}, computed tomography \cite{ref20, ref16, ref17}, and so on; see \cite{ref45} and the references therein. At the $k$-th iteration of the Kaczmarz method, the approximate solution ${\bf x}_{k+1}$ is obtained from orthogonally projecting the current iteration vector ${\bf x}_{k}$ onto the hyperplane  $H_{i_k}=\big\{ {\bf x}|A^{\left( i_k \right)}{\bf x}={\bf b}^{\left( i_k \right)} \big\}$:
\begin{equation}
	{\bf x}_{k+1}={\bf x}_k+\frac{{\bf b}^{\left ( i_k \right )}-A^{\left ( i_k \right )}{\bf x}_k}{\left \| A^{\left ( i_k \right )} \right \|_{2}^{2}}\left ( A^{\left ( i_k \right )} \right )^{H},\quad i_k = (k~mod~m) + 1.
\end{equation}

However, if the angle between the hyperplane related to two consecutive iterations is small, the convergence rate of Kaczmarz method may be slow. In order to partially overcome this difficulty, Strohmer and Vershynin \cite{ref2} proposed a randomized version and proved its expected exponential rate of convergence, which is called the randomized Kaczmarz (RK) method. However, RK is difficult to utilize some real-time information for selecting working rows in each iteration. Another disadvantage is that the probability of each row being selected is directly proportional to the square of its Euclidean norm. If the Euclidean norms of a few rows of the coefficient matrix are much larger than those of the others, the RK method will converge very slow or even diverge. In order to improve the convergence of the RK method, many variations have been proposed, such as the sampling Kaczmarz-Motzkin method \cite{ref22}, the accelerated randomized Kaczmarz method \cite{ref23}, the weighted randomized Kaczmarz method \cite{ref24}, and so on \cite{ref29,ref25, ref43, ref27, ref28, ref42}. Some randomized block Kaczmarz methods were introduced in \cite{ref38, ref37, ref35, ref36, ref31, ref30, ref32}.

In \cite{ref3}, Bai and Wu proposed a greedy randomized Kaczmarz (GRK) method for large consistent linear systems. Theoretical analysis and numerical experiments demonstrate that GRK works better than RK. However, in each iteration of the GRK method, one has to compute the residual vector corresponding to the augmented linear system, and construct a new index set, which suffer from a large amount of workload, especially for big data problems. In order to reduce the computational overhead, Jiang, Wu and Jiang \cite{ref4} proposed a semi-randomized Kaczmarz method with simple random sampling. This method only uses a small portion of rows of the matrix $A$, and there is no need to calculate probabilities nor constructing index set for choosing working rows.

On the other hand, if the linear system \eqref{1.1} is inconsistent, the Kaczmarz-type methods can not converge to the least-norm least-squares solution ${\bf x}_{\star}=A^{\dagger}{\bf b}$ \cite{ref6}. Inspired by the works in \cite{ref39, ref40}, Zouzias and Freris \cite{ref5} introduced a randomized extended Kaczmarz (REK) method, and proved that the REK method converges to the least-norm least-squares solution. Indeed, the REK method can be understood as using the randomized Kaczmarz method to solve two linear systems $A^H{\bf z}=0$ and $A{\bf x}={\bf b}-{\bf z}$. Based on the REK method, Bai and Wu \cite{ref13} proposed a partially random extended Kaczmarz (PREK) method. However, both the REK method and the PREK method will converge slowly if one of the two linear systems does so. A number of randomized block Kaczmarz methods for solving inconsistent linear systems were investigated in \cite{ref33, ref34}.

Recently, Bai and Wu proposed a greedy randomized augmented Kaczmarz (GRAK) method for solving  large sparse inconsistent linear systems \cite{ref6}. In essence, the GRAK method first equivalently transform the inconsistent linear system \eqref{1.1} into a consistent augmented linear system, and then apply the GRK method to solve this augmented linear system. The relaxed version of the GRAK method was proposed in \cite{ref12}. However, similar to the GRK method, it is required to calculate the residual vector with respect to the augmented linear system and construct two index sets and in each iteration. This is unfavorable for extremely large-scale problems.

In this paper, we focus on solving large (dense or sparse) inconsistent linear systems, and try to improve the performance of the GRAK method proposed in \cite{ref6}. The contributions of this work are as follows. First, we modify the iteration scheme for updating the approximate solution, and propose an accelerated greedy randomized augmented Kaczmarz method. We prove that it converges faster than the GRAK method if $\left\| A \right\| _{F}^{2}\gg 1$. Second, we propose a semi-randomized augmented Kaczmarz method with simple random sampling for large inconsistent linear systems. In this method, we are free of computing the residual vector corresponding to the augmented linear system, and there is no need to construct index sets. So the proposed method is much cheaper than the GRAK method per iteration. The convergence of the method is established. Third, as far as we know, there are no practical stopping criteria for Kaczmarz-type methods. The third contribution of this work is to introduce a practical stopping criterion for Kaczmarz-type methods. By using this strategy, there is no need to know the exact solution {\it a prior} nor compute the residual vector.

The organization of this paper is as follows. In Section 2, we briefly review the randomized extended Kaczmarz method and the greedy randomized augmented Kaczmarz method for large inconsistent linear systems. In Section 3, we propose an accelerated greedy randomized augmented Kaczmarz method and an accelerated greedy randomized augmented Kaczmarz method with simple random sampling. Convergence results of the two methods are established. In Section 4, we present a practical stopping criterion for Kaczmarz-type methods, and show its rationality theoretically. In Section 5, we perform comprehensive numerical experiments on both synthetic and real-world data sets, to show the efficiency of the proposed methods as well as the effectiveness of our stopping criterion.
Some concluding remarks are given in Section 6.

Let us introduce some notations. For a given matrix $A\in \mathbb{C} ^{m\times n}$, we denote by $A^H$, $\left\| A \right\| _2$, $A^{\dagger}$, $A^{\left( i \right)}$ and $A_{\left( j \right)}$ its conjugate transpose, Euclidean norm, Moore-Penrose inverse, the $i$th row and the $j$th column of $A$, respectively. Let $\mathcal{R} \left( A \right)$ and $\mathcal{R} \left( A \right) ^{\bot}$ be the range space of $A$ and the orthogonal complement subspace of $\mathcal{R} \left( A \right)$. Let ${\bf b}_{\mathcal{R} \left( A \right)}$ and ${\bf b}_{\mathcal{R} \left( A \right) ^{\bot}}$ be the orthogonal projection of the vector ${\bf b}$ onto $\mathcal{R} \left( A \right)$ and $\mathcal{R} \left( A \right) ^{\bot}$, respectively. Denote by $\lambda_{\min}(A^HA)$ the minimal nonzero eigenvalue of the matrix $A^HA$, and by ${\bf e_i}$ the $i$th column of the identity matrix $I$. Let $\mathbb{E} \left[ \cdot \right]$ be the full expectation, and let $\mathbb{E} _k\left[ \cdot \right]$ be the conditional expectation on the first $k$ iterations \cite{ref3}. That is,
    \begin{equation*}
    	\mathbb{E} _k\left[ \cdot \right] =\mathbb{E} \left[ \cdot |i_0,i_1,\ldots ,i_{k-1} \right],\quad  k=1,2,\ldots
    \end{equation*}
In the light of the properties of full expectation and conditional expectation, we have that $\mathbb{E} \left[ \mathbb{E} _k\left[ \cdot \right] \right] =\mathbb{E} \left[ \cdot \right]$.

\section{The randomized extended Kaczmarz method and the greedy randomized augmented Kaczmarz method}\label{Sec2}
\setcounter{equation}{0}

In this section, we briefly introduce the randomized extended Kaczmarz (REK) method and the greedy randomized augmented Kaczmarz (GRAK) method for large-scale inconsistent linear systems. The REK method is a combination of the randomized orthogonal projection algorithm and the randomized Kaczmarz method \cite{ref5}. The main idea is to efficiently reduce the norm of the ``noisy" of ${\bf b}$, i.e., ${\bf b}_{\mathcal{R} \left( A \right) ^{\bot}}$, by using the randomized orthogonal projection algorithm, and then apply the randomized Kaczmarz method on a new linear system whose right-hand side is now (arbitrarily) close to the column space of $A$. The framework of the REK method is as follows; for more details, refer to \cite{ref5}.

\begin{algorithm}\label{alg2}
	{\bf The randomized extended Kaczmarz~(REK) method}~{\rm\cite{ref5}}\\
	{{\bf Input:} Given $A$, $\bf b$, the convergence tolerance $\epsilon>0$, and let ${\bf x}_0=0$ and ${\bf z}_0={\bf b}$.}\\
	{{\bf Output:} The approximate solution $\widetilde{\bf x}$.}\\
	{\bf 1}. for $k=0,1,2,\ldots,$ do\\
	{\bf 2}. Pick $i_k\in \left\{ 1,2,\ldots,m \right\}$ with probability $q_i=\left\| A^{\left( i \right)} \right\| _{2}^{2}/\left\| A \right\| _{F}^{2},~ i\in \left\{ 1,2,\ldots,m \right\}$.\\
	{\bf 3}. Pick $j_k\in \left\{ 1,2,\ldots,n \right\}$ with probability $p_j=\left\| A^{\left( j \right)} \right\| _{2}^{2}/\left\| A \right\| _{F}^{2},~ j\in \left\{ 1,2,\ldots,n \right\}$.\\
	{\bf 4}. Set ${\bf z}_{k+1}={\bf z}_k-\frac{A_{\left( j_k \right)}^{H}{\bf z}_k}{\left\| A_{\left( j_k \right)} \right\| _{2}^{2}}A_{\left( j_k \right)}$.\\
	{\bf 5}. Set ${\bf x}_{k+1}={\bf x}_k-\frac{{\bf b}^{\left( i_k \right)}-{\bf z}_{k}^{\left( i_k \right)}-A^{\left( i_k \right)}{\bf x}_k}{\left\| A_{\left( j_k \right)} \right\| _{2}^{2}}\left( A^{\left( i_k \right)} \right)^H$.\\
	{\bf 6}. Check every $8\cdot \min \left( m, n \right)$ iterations and terminate if\\
	\[\frac{\left\| A{\bf x}_k-\left( {\bf b}-{\bf z}_k \right) \right\| _2}{\left\| A \right\| _F\left\| {\bf x}_k \right\| _2}\leqslant \varepsilon\,\,~~ and~~ \,\,\frac{\left\| A^H{\bf z}_k \right\| _2}{\left\| A \right\| _{F}^{2}\left\| {\bf x}_k \right\| _2}\leqslant \varepsilon.\]
	{\bf 7}. end for
\end{algorithm}

In \cite{ref6}, Bai and Wu proposed a greedy randomized augmented Kaczmarz (GRAK) method.
The main idea is to rewrite the large-scale inconsistent linear system $A{\bf x}={\bf b}$ into a consistent augmented linear system
\begin{equation}\label{2.3}
	\tilde{A}\tilde{{\bf x}}=\tilde{{\bf b}},
\end{equation}
where
\begin{equation}\label{eq22}
	\tilde{A}=\left( \begin{matrix}
		I&    A\\
		A^H&  0\\
		\end{matrix}\right) ,\quad \tilde{{\bf x}}=\left( \begin{array}{c}
		{\bf z}\\
		{\bf x}\\
		\end{array} \right),\quad and\quad\tilde{{\bf b}}=\left[ \begin{array}{c}
		{\bf b}\\
		{\bf 0}\\
		\end{array} \right],
\end{equation}
and then apply the greedy randomized augmented Kaczmarz (GRK) method \cite{ref3} to the augmented linear system \eqref{2.3}.

More precise, define the two index sets
\begin{equation}\label{2.5}
	\varOmega _{k}^{R}=\left\{ i_k\in \left\{ 1,2,\ldots ,m \right\} :\frac{\left| {\bf b}^{\left( i_k \right)}-{\bf z}_{k}^{\left( i_k \right)}-A^{\left( i_k \right)}{\bf x}_k \right|^2}{1+\left\| A^{\left( i_k \right)} \right\| _{2}^{2}}\geqslant \epsilon _k\left( \left\| {\bf b}-{\bf z}_k-A{\bf x}_k \right\| _{2}^{2}+\left\| A^H{\bf z}_k \right\| _{2}^{2} \right) \right\},
\end{equation}
and
\begin{equation}\label{2.6}
	\varOmega _{k}^{C}=\left\{ j_k\in \left\{ 1,2,\ldots ,n \right\}:\frac{\left| A_{\left( j_k \right)}^{H}{\bf z}_k \right|^2}{\left\| A_{\left( j_k \right)} \right\| _{2}^{2}}\geqslant \varepsilon _k\left( \left\| {\bf b}-{\bf z}_k-A{\bf x}_k \right\| _{2}^{2}+\left\| A^H{\bf z}_k \right\| _{2}^{2} \right) \right\},
\end{equation}
where
\begin{equation}\label{2.7}
	\epsilon _k=\max \left\{ \epsilon _{k}^{R},\epsilon _{k}^{C} \right\},
\end{equation}
with
\begin{equation}\label{2.8}
	\epsilon _{k}^{R}=\frac{1}{2}\left( \frac{1}{\left\| {\bf b}-{\bf z}_k-A{\bf x}_k \right\| _{2}^{2}+\left\| A^H{\bf z}_k \right\| _{2}^{2}}\underset{1\leqslant i_k\leqslant m}{\max}\left\{ \frac{\left| {\bf b}^{\left( i_k \right)}-{\bf z}_{k}^{\left( i_k \right)}-A^{\left( i_k \right)}{\bf x}_k \right|^2}{1+\left\| A^{\left( i_k \right)} \right\| _{2}^{2}} \right\} +\frac{1}{m+2\left\| A \right\| _{F}^{2}} \right)
\end{equation}
and
\begin{equation}\label{2.9}
	\epsilon _{k}^{C}=\frac{1}{2}\left( \frac{1}{\left\| {\bf b}-{\bf z}_k-A{\bf x}_k \right\| _{2}^{2}+\left\| A^H{\bf z}_k \right\| _{2}^{2}}\underset{1\leqslant j_k\leqslant n}{\max}\left\{ \frac{\left| A_{\left( j_k \right)}^{H}{\bf z}_k \right|^2}{\left\| A^{\left( j_k \right)} \right\| _{2}^{2}}\right\} +\frac{1}{m+2\left\| A \right\| _{F}^{2}} \right).
\end{equation}
Let
\begin{equation}\label{2.10}
	\tilde{{\bf r}}_{k}^{\left( i \right)}=\begin{cases}
		{\bf b}^{\left( i \right)}-{\bf z}_{k}^{\left( i \right)}-A^{\left( i \right)}{\bf x}_k&   if\,\, i\in \varOmega _{k}^{R},\\
		0&   otherwise,\\
	\end{cases}
\end{equation}
and
\begin{equation}\label{2.11}
	\tilde{{\bf s}}_{k}^{\left( j \right)}=\begin{cases}
		-A_{\left( j \right)}^{H}{\bf z}_k&   if\,\, j\in \varOmega _{k}^{C},\\
		0&   otherwise.\\
	\end{cases}
\end{equation}
The greedy randomized augmented Kaczmarz method is given as follows; for more details, refer to \cite{ref6}.
\begin{algorithm}\label{alg3}
	{\bf The greedy randomized augmented Kaczmarz ~(GRAK)~method}~{\rm\cite{ref6}}\\
	{{\bf Input:} $A$, $\bf b$, $l$, ${\bf x}_0$, and ${\bf z}_0$.}\\
	{{\bf Output:} ${\bf x}_l$ and ${\bf z}_l$.}\\
	{\bf 1.} for $k=0, 1, \ldots, l-1$ do\\
	{\bf 2.} Compute the tolerances $\epsilon _{k}^{R}$, $\epsilon _{k}^{C}$ and $\epsilon _k$ by \eqref{2.7}--\eqref{2.8}, and identify the index sets $\varOmega _{k}^{R}$ and $\varOmega _{k}^{C}$ by the criteria \eqref{2.5} and \eqref{2.6}, and determine the projected residuals $\tilde{\bf r}_k$ and $\tilde{\bf s}_k$ by the definitions \eqref{2.10} and \eqref{2.11}, respectively.\\
	{\bf 3.} Select an index $t_k$ satisfying $1\leqslant t_k\leqslant m+n$, with the probability
	\[Pr\left( index=t_k \right)=\begin{cases}
		\frac{\left| \tilde{{\bf r}}_{k}^{\left( t_k \right)} \right|^2}{\left\| \tilde{{\bf r}}_k \right\| _{2}^{2}+\left\| \tilde{{\bf s}}_k \right\| _{2}^{2}}&   if\,\, 1\leqslant t_k\leqslant m,\\
		\frac{\left| \tilde{{\bf s}}_{k}^{\left( t_k-m \right)} \right|^2}{\left\| \tilde{{\bf r}}_k \right\| _{2}^{2}+\left\| \tilde{{\bf s}}_k \right\| _{2}^{2}}&   if\,\, m+1\leqslant t_k\leqslant m+n.\\
	\end{cases}\]
	{\bf 4.} If $1\leqslant t_k\leqslant m$, then set $i_k=t_k$ and compute
	\[{\bf z}_{k+1}={\bf z}_k+\frac{\left( {\bf b}^{\left( i_k \right)}-{\bf z}_{k}^{\left( i_k \right)}-A^{\left( i_k \right)}{\bf x}_k \right)}{1+\left\| A^{\left( i_k \right)} \right\| _{2}^{2}}{\bf e}_{i_k}\]
	and\\
	\[{\bf x}_{k+1}={\bf x}_k+\frac{\left( {\bf b}^{\left( i_k \right)}-{\bf z}_{k}^{\left( i_k \right)}-A^{\left( i_k \right)}{\bf x}_k \right)}{1+\left\| A^{\left( i_k \right)} \right\| _{2}^{2}}\left( A^{\left( i_k \right)} \right) ^H\]
	else if~~ $m+1\leqslant t_k\leqslant m+n$, then set $j_k=t_k-m$, and compute
	\[{\bf z}_{k+1}={\bf z}_k-\frac{A_{\left( j_k \right)}^{H}{\bf z}_k}{\left\| A_{\left( j_k \right)} \right\| _{2}^{2}}A_{\left( j_k \right)} \quad and \quad {\bf x}_{k+1}={\bf x}_k\]
	{\bf 5.} endfor
\end{algorithm}

The following theorem indicates that the GRAK method is convergent in expectation to the least-norm least-square solution ${\bf x}_{\star}=A^{\dagger}{\bf b}$.
\begin{theorem}{\rm\cite{ref6}}\label{Thm2.2}
Starting from any initial vectors ${\bf x}_0\in range\left( A^H \right)$ and ${\bf z}_0\in \mathbb{C} ^m$, the iteration sequences $\left\{ {\bf x}_k \right\} _{k=0}^{\infty}$ and $\left\{ {\bf z}_k \right\} _{k=0}^{\infty}$, generated by the GRAK method, converge in expectation to the least-norm least-squares solution ${\bf x}_{\star}=A^{\dagger}{\bf b}$ of the linear system \eqref{1.1} and the orthogonal projection vector ${\bf z}_{\star}={\bf b}_{\mathcal{R} \left( A \right) ^{\bot}}$, respectively. Moreover, the global solution error in expectation with respect to both iteration sequences $\left\{ {\bf x}_k \right\} _{k=0}^{\infty}$ and $\left\{ {\bf z}_k \right\} _{k=0}^{\infty}$ obeys
	\begin{equation*}
		\mathbb{E} \left( \left\| {\bf x}_1-{\bf x}_{\star} \right\| _{2}^{2}+\left\| {\bf z}_1-{\bf z}_{\star} \right\| _{2}^{2} \right) \leqslant \zeta \left( \left\| {\bf x}_0-{\bf x}_{\star} \right\| _{2}^{2}+\left\| {\bf z}_0-{\bf z}_{\star} \right\| _{2}^{2} \right)
	\end{equation*}
	and
	\begin{equation}\label{eqn2.1}
		\mathbb{E} _k\left( \left\| {\bf x}_{k+1}-{\bf x}_{\star} \right\| _{2}^{2}+\left\| {\bf z}_{k+1}-{\bf z}_{\star} \right\| _{2}^{2} \right) \leqslant \beta \left( \left\| {\bf x}_k-{\bf x}_{\star} \right\| _{2}^{2}+\left\| {\bf z}_k-{\bf z}_{\star} \right\| _{2}^{2} \right), k=1,2\ldots
	\end{equation}
	where
	\begin{equation*}
		\zeta =1-\frac{\eta}{2\left\| A \right\| _{F}^{2}+m}
	\end{equation*}
	and
	\begin{equation}\label{2.15}
		\beta =1-\frac{1}{2}\left( \frac{2\left\| A \right\| _{F}^{2}+m}{\gamma}+1 \right) \frac{\eta}{2\left\| A \right\| _{F}^{2}+m}
	\end{equation}
	with
	\begin{equation}\label{2.16}
		\eta =\min \left\{ 1, \left( \sqrt{\lambda _{\min}\left( A^HA \right) +\frac{1}{4}}-\frac{1}{2} \right) ^2 \right\}
	\end{equation}
	and
	\begin{equation}\label{2.17}
		\gamma =2\left\| A \right\| _{F}^{2}+m-\min \left\{ 1+\underset{1\leqslant i\leqslant m}{\min}\left\| A^{\left( i \right)} \right\| _{2}^{2},\underset{1\leqslant j\leqslant n}{\min}\left\| A_{\left( j \right)} \right\| _{2}^{2} \right\}.
	\end{equation}
	As a result, it holds that
	\begin{equation*}
		\mathbb{E} \left( \left\| {\bf x}_k-{\bf x}_{\star} \right\| _{2}^{2}+\left\| {\bf z}_k-{\bf z}_{\star} \right\| _{2}^{2} \right) \leqslant \beta ^{k-1}\zeta\cdot \left( \left\| {\bf x}_0-{\bf x}_{\star} \right\| _{2}^{2}+\left\| {\bf z}_0-{\bf z}_{\star} \right\| _ {2}^{2} \right)
	\end{equation*}
\end{theorem}

Theorem \ref{Thm2.2} shows that the convergence rate of GRAK strongly depends on the minimum singular value and the Frobenius norm, the number of rows, and the size of the Euclidean norms of all rows and columns of the matrix $A$.
%

\section{The Proposed Methods}\label{Sec3}
\setcounter{equation}{0}
In each iteration of the GRAK method, one has to compute two residual vectors $\tilde{{\bf r}}_{k}$, $\tilde{{\bf s}}_{k}$ and the probabilities $PR(\cdot)$, to construct two new index sets $\epsilon _{k}^{R}$, $\epsilon _{k}^{C}$. Thus, one has to scan all the $m+n$ rows of the coefficient matrix $\tilde{A}$, which is very time consuming as the matrix $A$ is extremely large, and it is contrary to the original motivation of the Kaczmarz-type methods.

So as to deal with these problems, in this section, we first propose an accelerated greedy randomized augmented Kaczmarz method. Convergence analysis shows that, if $\|A\|_F^2\gg 1$, the accelerated method often converges faster than the GRAK method. To further reduce the computational overhead of the accelerated method, we then propose an accelerated greedy randomized augmented Kaczmarz method with simple sampling. The convergence of the method is also established.

\subsection{An accelerated greedy randomized augmented Kaczmarz method}
It is seen from Algorithm \ref{alg3} that, as $m+1\leqslant t_k\leqslant m+n$, the GRAK method only updates the vector ${\bf z}_k$ rather than the approximate solution ${\bf x}_k$, which may slow down the rate of convergence of this algorithm. Indeed, just like the REK method, GRAK may converge slowly if one of the iteration sequences $\left\{ {\bf x}_k \right\} _{k=0}^{\infty}$ or $\left\{ {\bf z}_k \right\} _{k=0}^{\infty}$ does so. Thus, it is necessary to give new insight into the GRAK method and improve its numerical performance.

The idea is that, when $m+1\leqslant t_k\leqslant m+n$, after computing ${\bf z}_k$, we select the working rows with probability $\| A^{\left( i_k \right)}\| _{2}^{2}/\| A \| _{F}^{2}$ to further update ${\bf x}_k$ for the linear equation $A{\bf x}={\bf b}-{\bf z}_{k+1}$.
Motivated by the semi-randomized Kaczmarz method proposed in \cite{ref4}, to refrain from evaluating probabilities and free of constructing two index sets in each iteration, we propose an accelerated greedy randomized augmented Kaczmarz (AGRAK) method to solve the augmented linear system \eqref{2.3}.
\begin{algorithm}\label{alg5}
	{\bf An accelerated greedy randomized augmented Kaczmarz (AGRAK) method}\\
	{{\bf Input:} $A$, $\bf b$, $l$, ${\bf x}_0={\bf 0}$, and ${\bf z}_0={\bf b}$.}\\
	{{\bf Output:} ${\bf x}_l$ and ${\bf z}_l$.}\\
		{\bf 1.} for $k=0, 1, \ldots, l-1$ do\\
		{\bf 2.} Select $t_k\in \left\{ 1,2,\dots ,m+n \right\}$ according to\[t_k=\mathop{\arg\max} \left\{ \underset{1\leqslant i\leqslant m}{\max}\left\{ \frac{\left| {\bf b}^{\left( i \right)}-{\bf z}_{k}^{\left( i \right)}-A^{\left( i \right)}{\bf x}_k \right|^2}{1+\left\| A^{\left( i \right)} \right\| _{2}^{2}}\right\}, \underset{m+1\leqslant i\leqslant m+n}{\max}\left\{ \frac{\left| A_{\left( i-m\right)}^{H}{\bf z}_k \right|^2}{\left\| A_{\left( i-m \right)} \right\| _{2}^{2}} \right\} \right\}.\]
		{\bf 3.} If $1\leqslant t_k\leqslant m$, then set $i_k=t_k$ and compute
		\[{\bf z}_{k+1}={\bf z}_k+\frac{\left( {\bf b}^{\left( i_k \right)}-{\bf z}_{k}^{\left( i_k \right)}-A^{\left( i_k \right)}{\bf x}_k \right)}{1+\left\| A^{\left( i_k \right)} \right\| _{2}^{2}}{\bf e}_{i_k},\]
		and\\
		\[{\bf x}_{k+1}={\bf x}_k+\frac{\left( {\bf b}^{\left( i_k \right)}-{\bf z}_{k}^{\left( i_k \right)}-A^{\left( i_k \right)}{\bf x}_k \right)}{1+\left\| A^{\left( i_k \right)} \right\| _{2}^{2}}\left( A^{\left( i_k \right)} \right) ^H .\]
		Else if $m+1\leqslant t_k\leqslant m+n$, then set $j_k=t_k-m$,
		\[{\bf z}_{k+1}={\bf z}_k-\frac{A_{\left( j_k \right)}^{H}{\bf z}_k}{\left\| A_{\left( j_k \right)} \right\| _{2}^{2}}A_{\left( j_k \right)},\]
		and select an index $i_k\in \left\{ 1,2,\dots ,m \right\}$ with probability $P_r\left( row=i_k \right) =\frac{\left\| A^{\left( i_k \right)} \right\| _{2}^{2}}{\left\| A \right\| _{F}^{2}}$, and compute
		\[{\bf x}_{k+1}={\bf x}_k+\frac{\left( {\bf b}^{\left( i_k \right)}-{\bf z}_{k+1}^{\left( i_k \right)}-A^{\left( i_k \right)}{\bf x}_k \right)}{\left\| A^{\left( i_k \right)} \right\| _{2}^{2}}\left( A^{\left( i_k \right)} \right) ^H.\]
		{\bf 4.} endfor
\end{algorithm}

Next, we consider the convergence of the AGRAK method. The following lemma is needed.
\begin{Lemma}\label{Lem4.1}
	Let $a_1,a_2,b_1,b_2\in \mathbb{R} ^+$, and $a_1<a_2$. Then there exists $\upsilon \in \left( a_1,a_2 \right)$, such that $a_1b_1+a_2b_2=\upsilon \left( b_1+b_2 \right)$.
\end{Lemma}
\begin{proof}
	Let $\upsilon =\frac{a_1b_1+a_2b_2}{b_1+b_2}$, then \[\upsilon =\frac{a_1b_1+a_2b_2}{b_1+b_2}<\frac{a_2b_1+a_2b_2}{b_1+b_2}=\frac{a_2\left( b_1+b_2 \right)}{b_1+b_2}=a_2,\]
	and\[\upsilon =\frac{a_1b_1+a_2b_2}{b_1+b_2}>\frac{a_1b_1+a_1b_2}{b_1+b_2}=\frac{a_1\left( b_1+b_2 \right)}{b_1+b_2}=a_1,\]
which satisfies $\upsilon \left( b_1+b_2 \right) =a_1b_1+a_2b_2$.
\end{proof}


 First, notice that if $1\leqslant t_k\leqslant m$, Algorithm \ref{alg5} applies the semi-randomized Kaczmarz method \cite{ref4} to solve the augmented linear system \eqref{2.3}. Denote by
 \begin{equation*}
       \tilde{{\bf x}}_{k}=\left[ \begin{array}{c}
		{\bf z}_{k}\\
		{\bf x}_{k}\\
		\end{array} \right]\quad and\quad\tilde{{\bf x}}_{\star}=\left[ \begin{array}{c}
		{\bf z}_{\star}\\
		{\bf x}_{\star}\\
		\end{array} \right],
\end{equation*}
 from \cite[Theorem 3.2]{ref4}, we have that
\begin{equation}\label{4.7}
	\mathbb{E} _k\left\| \tilde{{\bf x}}_{k+1}-\tilde{{\bf x}}_{\star} \right\| _{2}^{2}\leqslant \left( 1-\frac{\lambda _{\min}\left( \tilde{A}^H\tilde{A} \right)}{\left\| \tilde{A} \right\| _{F}^{2}-\underset{1\leqslant i\leqslant m+n}{\min}\left\{ \left\| \tilde{A}^{\left( i \right)} \right\| _{2}^{2} \right\}} \right) \left\| \tilde{{\bf x}}_k-\tilde{{\bf x}}_{\star} \right\| _{2}^{2}.
\end{equation}
Thanks to the structure of $\tilde{A}$, $\tilde{{\bf x}}_{k+1}$ and $\tilde{{\bf x}}_{\star}$; see \eqref{eq22}, one can rewrite \eqref{4.7} as
\begin{equation}\label{4.8}
	\mathbb{E} _k\left( \left\| {\bf x}_{k+1}-{\bf x}_{\star} \right\| _{2}^{2}+\left\| {\bf z}_{k+1}-{\bf z}_{\star} \right\| _{2}^{2} \right) \leqslant \left( 1-\frac{\eta}{\gamma} \right) \left( \left\| {\bf x}_k-{\bf x}_{\star} \right\| _{2}^{2}+\left\| {\bf z}_k-{\bf z}_{\star} \right\| _{2}^{2} \right),
\end{equation}
where $\eta$ and $\gamma$ be defined in \eqref{2.16} and \eqref{2.17}.

Take full expectation on the both sides of \eqref{4.8}, and set
\begin{equation}\label{33}
\tilde{\beta}=1-\frac{\eta}{\gamma},
\end{equation}
we arrive at
\begin{equation*}
	\mathbb{E} \left( \left\| {\bf x}_{k+1}-{\bf x}_{\star} \right\| _{2}^{2}+\left\| {\bf z}_{k+1}-{\bf z}_{\star} \right\| _{2}^{2} \right) \leqslant \tilde{\beta}\mathbb{E} \left( \left\| {\bf x}_k-{\bf x}_{\star} \right\| _{2}^{2}+\left\| {\bf z}_k-{\bf z}_{\star} \right\| _{2}^{2} \right), \quad 1\leqslant t_k\leqslant m.
\end{equation*}

Second, we consider the convergence of ${\bf z}_k$ as $m+1\leqslant t_k\leqslant m+n$. Notice that if $m+1\leqslant t_k\leqslant m+n$, then
\begin{equation*}
	\max \left\{ \underset{1\leqslant i\leqslant m}{\max}\left\{ \frac{\left| {\bf b}^{\left( i \right)}-{\bf z}_{k}^{\left( i \right)}-A^{\left( i \right)}{\bf x}_k \right|^2}{1+\left\| A^{\left( i \right)} \right\| _{2}^{2}}\right\}, \underset{m+1\leqslant i\leqslant m+n}{\max}\left\{ \frac{\left| A_{\left( i-m\right)}^{H}{\bf z}_k \right|^2}{\left\| A_{\left( i-m \right)} \right\| _{2}^{2}} \right\} \right\} = \underset{m+1\leqslant i\leqslant m+n}{\max}\left\{ \frac{\left| A_{\left( i-m\right)}^{H}{\bf z}_k \right|^2}{\left\| A_{\left( i-m \right)} \right\| _{2}^{2}} \right\}.
\end{equation*}
Consequently,
\begin{equation*}
	j_k=t_k-m=\mathop{\arg\max}_{1\leqslant j\leqslant n}\left\{ \frac{\left| A_{\left( j \right)}^{H}{\bf z}_k \right|^2}{\left\| A_{\left( j \right)} \right\| _{2}^{2}} \right\},
\end{equation*}
 and
\begin{equation*}
	{\bf z}_{k+1}={\bf z}_k-\frac{A_{\left( j_k \right)}^{H}{\bf z}_k}{\left\| A_{\left( j_k \right)} \right\| _{2}^{2}}A_{\left( j_k \right)}, \quad  where\quad j_k=\mathop{\arg\max}_{1\leqslant j\leqslant n}\left\{ \frac{\left| A_{\left( j \right)}^{H}{\bf z}_k \right|^2}{\left\| A_{\left( j \right)} \right\| _{2}^{2}} \right\},
\end{equation*}
 which is nothing but applying the semi-randomized Kaczmarz method to solve the equation $A^H{\bf z}=0$. Therefore, by \cite[Theorem 3.2]{ref4},
\begin{equation}\label{4.11}
	\mathbb{E} _k\left\| {\bf z}_{k+1}-{\bf z}_{\star} \right\| _{2}^{2}\leqslant \left( 1-\frac{\lambda _{\min}\left( A^HA \right)}{\left\| A \right\| _{F}^{2}-\underset{1\leqslant j\leqslant n}{\min}\left\{ \left\| A_{\left( j \right)} \right\| _{2}^{2} \right\}} \right) \left\| {\bf z}_k-{\bf z}_{\star} \right\| _{2}^{2}, \quad m+1\leqslant t_k\leqslant m+n.
\end{equation}
Take full expectation on the both sides of \eqref{4.11}, and set
\begin{equation*}
	\delta =1-\frac{\lambda _{\min}\left( A^HA \right)}{\left\| A \right\| _{F}^{2}-\underset{1\leqslant j\leqslant n}{\min}\left\{ \left\| A_{\left( j \right)} \right\| _{2}^{2} \right\}},
\end{equation*}
we get
\begin{equation}\label{4.12}
	\mathbb{E} \left\| {\bf z}_{k+1}-{\bf z}_{\star} \right\| _{2}^{2}\leqslant \delta \mathbb{E} \left\| {\bf z}_k-{\bf z}_{\star} \right\| _{2}^{2}.
\end{equation}
Next, we  update ${\bf x}_k$ for the linear equation $A{\bf x}={\bf b}-{\bf z}_{k+1}$ by using the randomized Kaczmarz method. From \cite{ref44, ref2}, we obtain
\begin{equation}\label{4.13}
	\mathbb{E} _k\left\| {\bf x}_{k+1}-{\bf x}_{\star} \right\| _{2}^{2}\leqslant \left( 1-\frac{\lambda _{\min}\left( A^HA \right)}{\left\| A \right\| _{F}^{2}} \right) \left\| {\bf x}_k-{\bf x}_{\star} \right\| _{2}^{2}+\frac{\left\| {\bf z}_{k+1}-{\bf z}_{\star} \right\| _{2}^{2}}{\left\| A \right\| _{F}^{2}}, \quad m+1\leqslant t_k\leqslant m+n.
\end{equation}
Take full expectation on the both sides of \eqref{4.13}, and set
\begin{equation}\label{alp}
	\alpha =1-\frac{\lambda _{\min}\left( A^HA \right)}{\left\| A \right\| _{F}^{2}},
\end{equation}
we get
\begin{align}\label{4.14}
\mathbb{E} \left\| {\bf x}_{k+1}-{\bf x}_{\star} \right\| _{2}^{2}
		&\leqslant \alpha \mathbb{E} \left\| {\bf x}_k-{\bf x}_{\star} \right\| _{2}^{2}+\frac{\mathbb{E} \left\| {\bf z}_{k+1}-{\bf z}_{\star} \right\| _{2}^{2}}{\left\| A \right\| _{F}^{2}}\nonumber\\
		&\leqslant \alpha \mathbb{E} \left\| {\bf x}_k-{\bf x}_{\star} \right\| _{2}^{2}+\frac{\delta}{\left\| A \right\| _{F}^{2}}\mathbb{E} \left\| {\bf z}_k-{\bf z}_{\star} \right\| _{2}^{2},
\end{align}
where the second inequality is from \eqref{4.12}.
Suppose that $\left\| A \right\| _F^2\gg 1$\footnote{We mention that this assumption is very weak in practice, since Kaczmarz-type methods often applies to very large (dense) linear systems.}, such that
\begin{equation*}
	\frac{\delta}{\left\| A \right\| _{F}^{2}}\mathbb{E} \left\| {\bf z}_k-{\bf z}_{\star} \right\| _{2}^{2}\ll \alpha \mathbb{E} \left\| {\bf x}_k-{\bf x}_{\star} \right\| _{2}^{2}+\delta \mathbb{E} \left\| {\bf z}_k-{\bf z}_{\star} \right\| _{2}^{2},
\end{equation*}
or in other words,
\begin{equation}\label{4.15}
	\frac{\delta}{\left\| A \right\| _{F}^{2}}\mathbb{E} \left\| {\bf z}_k-{\bf z}_{\star} \right\| _{2}^{2} =o \left(  \alpha \mathbb{E} \left\| {\bf x}_k-{\bf x}_{\star} \right\| _{2}^{2}+\delta \mathbb{E} \left\| {\bf z}_k-{\bf z}_{\star} \right\| _{2}^{2} \right).
\end{equation}
Combining \eqref{4.12}, \eqref{4.14} and \eqref{4.15}, we have from Lemma \ref{Lem4.1} that there exists  $\delta <\theta _k <\alpha$, such that
\begin{align}\label{4.24}
	\mathbb{E} \left( \left\| {\bf x}_{k+1}-{\bf x}_{\star} \right\| _{2}^{2}+\left\| {\bf z}_{k+1}-{\bf z}_{\star} \right\| _{2}^{2} \right) &\leqslant \alpha \mathbb{E} \left\| {\bf x}_k-{\bf x}_{\star} \right\| _{2}^{2}+\delta \mathbb{E} \left\| {\bf z}_k-{\bf z}_{\star} \right\| _{2}^{2}+\frac{\delta}{\left\| A \right\| _{F}^{2}}\mathbb{E} \left\| {\bf z}_k-{\bf z}_{\star} \right\| _{2}^{2} \nonumber\\
	&=\alpha \mathbb{E} \left\| {\bf x}_k-{\bf x}_{\star} \right\| _{2}^{2}+\delta \mathbb{E} \left\| {\bf z}_k-{\bf z}_{\star} \right\| _{2}^{2}+o \big(\alpha \mathbb{E} \left\| {\bf x}_k-{\bf x}_{\star} \right\| _{2}^{2}+\delta \mathbb{E} \left\| {\bf z}_k-{\bf z}_{\star} \right\| _{2}^{2} \big) \nonumber\\
	&\approx \alpha \mathbb{E} \left\| {\bf x}_k-{\bf x}_{\star} \right\| _{2}^{2}+\delta \mathbb{E} \left\| {\bf z}_k-{\bf z}_{\star} \right\| _{2}^{2} \nonumber\\
	&= \theta _k\mathbb{E} \left( \left\| {\bf x}_k-{\bf x}_{\star} \right\| _{2}^{2}+\left\| {\bf z}_k-{\bf z}_{\star} \right\| _{2}^{2} \right), \quad m+1\leqslant t_k\leqslant m+n,
\end{align}
where the high order term $o \big(  \alpha \mathbb{E} \left\| {\bf x}_k-{\bf x}_{\star} \right\| _{2}^{2}+\delta \mathbb{E} \left\| {\bf z}_k-{\bf z}_{\star} \right\| _{2}^{2} \big)$ is omitted.
In conclusion, we have the follow theorem for the convergence of the proposed AGRAK method.
\begin{theorem}\label{Thm4.2}
	Let $\left\| A \right\| _{F}^{2}\gg 1$ such that \eqref{4.15} is satisfied. Then under the above notations, the iteration sequences $\left\{ {\bf x}_k \right\} _{k=1}^{\infty}$ and $\left\{ {\bf z}_k \right\} _{k=1}^{\infty}$ generated by Algorithm \ref{alg5} converge in expectation to ${\bf x}_{\star}$ and ${\bf z}_{\star}$, with
	\begin{equation*}
		\begin{cases}
			\mathbb{E} \left( \left\| {\bf x}_{k+1}-{\bf x}_{\star} \right\| _{2}^{2}+\left\| {\bf z}_{k+1}-{\bf z}_{\star} \right\| _{2}^{2} \right) \leqslant \tilde{\beta}\mathbb{E} \left( \left\| {\bf x}_k-{\bf x}_{\star} \right\| _{2}^{2}+\left\| {\bf z}_k-{\bf z}_{\star} \right\| _{2}^{2} \right)&      if\,\, 1\leqslant t_k\leqslant m,\\
			\mathbb{E} \left( \left\| {\bf x}_{k+1}-{\bf x}_{\star} \right\| _{2}^{2}+\left\| {\bf z}_{k+1}-{\bf z}_{\star} \right\| _{2}^{2} \right) \lesssim\theta _k\mathbb{E} \left( \left\| {\bf x}_k-{\bf x}_{\star} \right\| _{2}^{2}+\left\| {\bf z}_k-{\bf z}_{\star} \right\| _{2}^{2} \right)&      if\,\, m+1\leqslant t_k\leqslant m+n,\\
		\end{cases}
	\end{equation*}
where ``$\lesssim$" means the high order term $o \big(  \alpha \mathbb{E} \left\| {\bf x}_k-{\bf x}_{\star} \right\| _{2}^{2}+\delta \mathbb{E} \left\| {\bf z}_k-{\bf z}_{\star} \right\| _{2}^{2} \big)$ is omitted.
\end{theorem}

Under the condition that $\|A\|_F^2\gg 1$, the following theorem indicates that Algorithm \ref{alg5} converges faster than the GRAK method; refer to \eqref{eqn2.1}.
\begin{theorem}
Let $\beta,\tilde{\beta}$ and $\theta _k$ be defined in \eqref{2.15}, \eqref{33} and \eqref{4.24}, respectively. If $\left\| A \right\| _{F}^{2}\gg 1$ such that \eqref{4.15} is satisfied, then
\begin{equation*}
0<\tilde{\beta}<\beta\quad{\rm and}\quad 0<\theta _k<\beta.
\end{equation*}
\end{theorem}

\begin{proof}
First, we prove that $\tilde{\beta}< \beta$. Recall that
\begin{equation*}
  \beta =1-\frac{1}{2}\left( \frac{\left\| \tilde{A} \right\| _{F}^{2}}{\left\| \tilde{A} \right\| _{F}^{2}-\underset{1\leqslant i\leqslant m+n}{\min}\left\{ \left\| \tilde{A}^{\left( i \right)} \right\| _{2}^{2} \right\}}+1 \right) \frac{\lambda _{\min}\left( \tilde{A}^H\tilde{A} \right)}{\left\| \tilde{A} \right\| _{F}^{2}},
\end{equation*}
and
\begin{align*}
  \tilde{\beta}&=1-\frac{\lambda _{\min}\left( \tilde{A}^H\tilde{A} \right)}{\left\| \tilde{A} \right\| _{F}^{2}-\underset{1\leqslant i\leqslant m+n}{\min}\left\{ \left\| \tilde{A}^{\left( i \right)} \right\| _{2}^{2} \right\}}=1-\frac{\left\| \tilde{A} \right\| _{F}^{2}}{\left\| \tilde{A} \right\| _{F}^{2}-\underset{1\leqslant i\leqslant m+n}{\min}\left\{ \left\| \tilde{A}^{\left( i \right)} \right\| _{2}^{2} \right\}}\cdot \frac{\lambda _{\min}\left( \tilde{A}^H\tilde{A} \right)}{\left\| \tilde{A} \right\| _{F}^{2}}.
\end{align*}
As there are no zero columns in $A$, we obtain
\begin{equation*}
   \frac{1}{2}\left( \frac{\left\| \tilde{A} \right\| _{F}^{2}}{\left\| \tilde{A} \right\| _{F}^{2}-\underset{1\leqslant i\leqslant m+n}{\min}\left\{ \left\| \tilde{A}^{\left( i \right)} \right\| _{2}^{2} \right\}}+1 \right) <\frac{\left\| \tilde{A} \right\| _{F}^{2}}{\left\| \tilde{A} \right\| _{F}^{2}-\underset{1\leqslant i\leqslant m+n}{\min}\left\{ \left\| \tilde{A}^{\left( i \right)} \right\| _{2}^{2} \right\}},
\end{equation*}
and $\tilde{\beta}<\beta$.

Second, we prove that $\theta _k <\beta$. Indeed, it is sufficient to show that $\alpha <\beta$. Let
\begin{equation*}
	\begin{cases}
		\eta =1<\lambda _{\min}\left( A^HA \right) ,&    if\,\, \lambda _{\min}\left( A^HA \right) \geqslant 2,\\
		\eta =\left( \sqrt{\lambda _{\min}\left( A^HA \right) +\frac{1}{4}}-\frac{1}{2} \right) ^2,&     if\,\, \lambda _{\min}\left( A^HA \right) <2.\\
	\end{cases}
\end{equation*}
Notice that $\eta <\lambda _{\min}\left( A^HA \right)$.
Let $\xi =2+\frac{m}{\left\| A \right\| _{F}^{2}}$, then
\begin{equation}\label{4.16}
	\frac{1}{2}\left( \frac{2\left\| A \right\| _{F}^{2}+m}{\gamma}+1 \right) \frac{\eta}{2\left\| A \right\| _{F}^{2}+m}=\frac{1}{2\xi}\left( \frac{2\left\| A \right\| _{F}^{2}+m}{\gamma}+1 \right) \frac{\eta}{\left\| A \right\| _{F}^{2}}.
\end{equation}

Next, we show that
$\frac{2\left\| A \right\| _{F}^{2}+m}{\gamma}+1<4.$
Otherwise, if
$\frac{2\left\| A \right\| _{F}^{2}+m}{\gamma}+1\geqslant4$,
then
\begin{equation*}
	3\left\| A \right\| _{F}^{2}\geqslant 3\cdot \min \left\{ 1+\underset{1\leqslant i\leqslant m}{\min}\left\{ \left\| A^{\left( i \right)} \right\| _{2}^{2} \right\} , \underset{1\leqslant j\leqslant n}{\min}\left\{ \left\| A_{\left( j \right)} \right\| _{2}^{2} \right\} \right\}\geq 4\left\| A \right\| _{F}^{2} +2m,
\end{equation*}
which is incorrect. As a result,
\begin{equation*}
	\frac{1}{2\xi}\left( \frac{2\left\| A \right\| _{F}^{2}+m}{\gamma}+1 \right) <1,
\end{equation*}
and
\begin{equation}\label{4.25}
	\frac{1}{2\xi}\left( \frac{2\left\| A \right\|_{F}^{2}+m}{\gamma}+1 \right) \frac{\eta}{\left\| A \right\|_{F}^{2}}<\frac{\lambda _{\min}\left( A^HA \right)}{\left\| A \right\| _{F}^{2}}.
\end{equation}
A combination of \eqref{4.16} and \eqref{4.25} yields
\begin{equation*}
	\frac{1}{2}\left( \frac{2\left\| A \right\| _{F}^{2}+m}{\gamma}+1 \right) \frac{\eta}{2\left\| A \right\| _{F}^{2}+m}<\frac{\lambda _{\min}\left( A^HA \right)}{\left\| A \right\| _{F}^{2}},
\end{equation*}
and it follows from \eqref{2.15} and \eqref{alp} that
$\theta _k <\alpha <\beta.$
\end{proof}


\subsection{An accelerated greedy randomized augmented Kaczmarz method with simple random sampling}

Although Algorithm \ref{alg5} can converge faster than GRAK, however, just like the GRAK method, one has to compute the residual vector $ A{\bf x}_k-( {\bf b}-{\bf z}_k )$ and $A^H{\bf z}_k$ in each iteration. That is, we have to access all the rows and columns of the data matrix $A$, which is contrary to the motivation of the Kaczmarz-type methods. Specifically, this is unfavorable when the coefficient matrix is large and dense, and it is interesting to accelerate Algorithm \ref{alg5} further.

To deal with this problem, we apply the simple random sampling strategy advocated in \cite{ref4} to solve \eqref{2.3}. The key idea is that, in each iteration, we first generate a subset by using simple random sampling, and then select working rows in this subset. Notice that the cost of generating a subset via simple random sampling is negligible compared with that of the overall algorithm.
The advantages are that we are free of accessing all the rows or columns of the data matrix, and there is no need to compute the residual vector associated with the augmented linear system. The algorithm is given as follows.

\begin{algorithm}\label{alg6}
	{\bf A semi-randomized augmented Kaczmarz method with simple random sampling}\\
	{{\bf Input:} $A$, $\bf b$, ${\bf x}_0$, ${\bf z}_0$, a parameter $0<\eta \ll 1$ and the maximal iteration number $l$.}\\
	{{\bf Output:} ${\bf x}_l$ and ${\bf z}_l$.}\\
		{\bf 1.} for $k=0, 1, \ldots, l-1$ do\\
		{\bf 2.} Generate an indicator set $\varOmega _k$ with $\lfloor\left( m+n \right) \eta\rceil$ rows by using the simple random sampling method, where $\lfloor\cdot\rceil$ means rounding down. Let $\varOmega _{k}^{1}=\varOmega _k\cap \left\{ 1,2,\dots ,m \right\}$ and $\varOmega _{k}^{2}=\varOmega _k\cap \left\{ m+1,m+2,\dots ,m+n \right\}$. Then, select $t_k$ according to
		\begin{equation*}
			t_k=\mathop{\arg\max} \left\{ \max_{j_k\in \varOmega _{k}^{1}} \left\{ \frac{\left| {\bf b}^{\left( j_k \right)}-{\bf z}_{k}^{\left( j_k \right)}-A^{\left( j_k \right)}{\bf x}_k \right|}{\sqrt{1+\left\| A^{\left( j_k \right)} \right\| _{2}^{2}}} \right\} ,\max_{j_k\in \varOmega _{k}^{2}} \left\{ \frac{\left| -A_{\left( j_k-m \right)}^{H}{\bf z}_k \right|}{\left\| A_{\left( j_k-m \right)}^{H} \right\| _2} \right\} \right\}
		\end{equation*}
		{\bf 3.} If $1\leqslant t_k\leqslant m$, then set $i_k=t_k$, and compute
		\[{\bf z}_{k+1}={\bf z}_k+\frac{\left( {\bf b}^{\left( i_k \right)}-{\bf z}_{k}^{\left( i_k \right)}-A^{\left( i_k \right)}{\bf x}_k \right)}{1+\left\| A^{\left( i_k \right)} \right\| _{2}^{2}}{\bf e}_{i_k},\]
		and
		\[{\bf x}_{k+1}={\bf x}_k+\frac{\left( {\bf b}^{\left( i_k \right)}-{\bf z}_{k}^{\left( i_k \right)}-A^{\left( i_k \right)}{\bf x}_k \right)}{1+\left\| A^{\left( i_k \right)} \right\| _{2}^{2}}\left( A^{\left( i_k \right)} \right) ^H.\]
		Else if $m+1\leqslant t_k\leqslant m+n$, then set $j_k=t_k-m$, and compute
		\begin{equation*}
			{\bf z}_{k+1}={\bf z}_k-\frac{A_{\left( j_k \right)}^{H}{\bf z}_k}{\left\| A_{\left( j_k \right)} \right\| _{2}^{2}}A_{\left( j_k \right)},
		\end{equation*}
		and select an index $i_k\in \left\{ 1,2,\dots ,m \right\}$ with probability $P_r\left( row=i_k \right) =\frac{\left\| A^{\left( i_k \right)} \right\| _{2}^{2}}{\left\| A \right\| _{F}^{2}}$, and compute
		\[{\bf x}_{k+1}={\bf x}_k+\frac{\left( {\bf b}^{\left( i_k \right)}-{\bf z}_{k+1}^{\left( i_k \right)}-A^{\left( i_k \right)}{\bf x}_k \right)}{\left\| A^{\left( i_k \right)} \right\| _{2}^{2}}\left( A^{\left( i_k \right)} \right) ^H.\]
		 {\bf 4.} end for
	\end{algorithm}

We are in a position to consider the convergence of the semi-randomized augmented Kaczmarz method with simple random sampling. First, if $1\leqslant t_k\leqslant m$, similar to the proof of Theorem 4.2 in \cite{ref4}, as the number of the rows of $A$ is sufficiently large, we have from the Chebyshev's (weak) law of large numbers that \cite{WLaw}, there are two scalars $0\leqslant \hat{\varepsilon}_k,\tilde{\varepsilon}_k\ll 1$ such that

\begin{equation*}
	\mathbb{E} _k\left\| \tilde{{\bf x}}_{k+1}-\tilde{{\bf x}}_{\star} \right\| _{2}^{2}\leqslant \left( 1-\frac{\left( 1-\hat{\varepsilon}_k \right) \lambda _{\min}\left( \tilde{A}^H\tilde{A} \right)}{\left( 1+\tilde{\varepsilon}_k \right) \left( \left\| \tilde{A} \right\| _{F}^{2}-\underset{1\leqslant i\leqslant m+n}{\min}\left\{ \left\| \tilde{A}^{\left( i \right)} \right\| _{2}^{2} \right\} \right)} \right) \left\| \tilde{{\bf x}}_k-\tilde{{\bf x}}_{\star} \right\| _{2}^{2}.
\end{equation*}
It can be rewritten as
\begin{equation}\label{4.19}
	\mathbb{E} _k\left( \left\| {\bf x}_{k+1}-{\bf x}_{\star} \right\| _{2}^{2}+\left\| {\bf z}_{k+1}-{\bf z}_{\star} \right\| _{2}^{2} \right)\leqslant \left( 1-\frac{\left( 1-\hat{\varepsilon}_k \right) \eta}{\left( 1+\tilde{\varepsilon}_k \right) \gamma} \right) \left( \left\| {\bf x}_k-{\bf x}_{\star} \right\| _{2}^{2}+\left\| {\bf z}_k-{\bf z}_{\star} \right\| _{2}^{2} \right),
\end{equation}
where $\eta$ and $\gamma$ are defined in \eqref{2.16} and \eqref{2.17}, respectively.
Take full expectation on the both sides of \eqref{4.19} and set
\begin{equation*}
	\hat{\beta}=1-\frac{\left( 1-\hat{\varepsilon}_k \right) \eta}{\left( 1+\tilde{\varepsilon}_k \right) \gamma},
\end{equation*}
we obtain
\begin{equation*}
	\mathbb{E} \left( \left\| {\bf x}_{k+1}-{\bf x}_{\star} \right\| _{2}^{2}+\left\| {\bf z}_{k+1}-{\bf z}_{\star} \right\| _{2}^{2} \right)\leqslant \hat{\beta} \mathbb{E} \left( \left\| {\bf x}_k-{\bf x}_{\star} \right\| _{2}^{2}+\left\| {\bf z}_k-{\bf z}_{\star} \right\| _{2}^{2} \right),\quad 1\leqslant t_k\leqslant m.
\end{equation*}

Second, if $m+1\leqslant t_k\leqslant m+n$, analogous to the analysis of Algorithm \ref{alg5}, we have that
\begin{equation*}
	\mathbb{E} \left\| {\bf z}_{k+1}-{\bf z}_{\star} \right\| _{2}^{2}\leqslant \tilde{\delta} \mathbb{E} \left\| {\bf z}_k-{\bf z}_{\star} \right\|_{2}^{2}.
\end{equation*}
As the number of the rows of $A$ is sufficiently large, we have from the Chebyshev's (weak) law of large numbers that \cite{WLaw}, there are two scalars $0\leqslant \check{\varepsilon}_{k},\mathring{\varepsilon}_{k}\ll 1$, such that
\begin{equation*}
	\mathbb{E} \left\| {\bf x}_{k+1}-{\bf x}_{\star} \right\| _{2}^{2}\leqslant \alpha \mathbb{E} \left\| {\bf x}_k-{\bf x}_{\star} \right\| _{2}^{2}+\frac{\tilde{\delta}}{\left\| A \right\| _{F}^{2}}\mathbb{E} \left\| {\bf z}_k-{\bf z}_{\star} \right\| _{2}^{2},
\end{equation*}
where $\alpha$ is defined in \eqref{alp}, and
\begin{equation*}
\tilde{\delta}=1-\frac{\left( 1-\check{\varepsilon}_{k} \right) \lambda _{\min}\left( A^HA \right)}{\left( 1+\mathring{\varepsilon}_{k} \right) \left( \left\| A \right\| _{F}^{2}-\underset{1\leqslant j\leqslant n}{\min}\left\{ \left\| A_{\left( j \right)} \right\| _{2}^{2} \right\} \right)}.
\end{equation*}

If $\left\| A \right\| _{F}^{2}\gg 1$ such that
\begin{equation*}
	\frac{\tilde{\delta}}{\left\| A \right\| _{F}^{2}}\mathbb{E} \left\| {\bf z}_k-{\bf z}_{\star} \right\| _{2}^{2}\ll \alpha \mathbb{E} \left\| {\bf x}_k-{\bf x}_{\star} \right\| _{2}^{2}+\tilde{\delta} \mathbb{E} \left\| {\bf z}_k-{\bf z}_{\star} \right\| _{2}^{2},
\end{equation*}
i.e.,
\begin{equation}\label{4.23}
	\frac{\tilde{\delta}}{\left\| A \right\| _{F}^{2}}\mathbb{E} \left\| {\bf z}_k-{\bf z}_{\star} \right\| _{2}^{2} =o \left(  \alpha \mathbb{E} \left\| {\bf x}_k-{\bf x}_{\star} \right\| _{2}^{2}+\tilde{\delta} \mathbb{E} \left\| {\bf z}_k-{\bf z}_{\star} \right\| _{2}^{2} \right).
\end{equation}
Then we have the following theorem on the convergence of Algorithm \ref{alg6}.

\begin{theorem}\label{4.3}
		Suppose that $\left\| A \right\| _{F}^{2}\gg 1$ such that \eqref{4.23} is satisfied. Then the iteration sequences $\left\{ {\bf x}_k \right\} _{k=1}^{\infty}$ and $\left\{ {\bf z}_k \right\} _{k=1}^{\infty}$ generated by Algorithm \ref{alg6}, converge in expectation to ${\bf x}_{\star}$ and ${\bf z}_{\star}$, with
	\begin{equation*}
		\begin{cases}
			\mathbb{E} \left( \left\| {\bf x}_{k+1}-{\bf x}_{\star} \right\| _{2}^{2}+\left\| {\bf z}_{k+1}-{\bf z}_{\star} \right\| _{2}^{2} \right)\leqslant \hat{\beta} \mathbb{E} \left( \left\| {\bf x}_k-{\bf x}_{\star} \right\| _{2}^{2}+\left\| {\bf z}_k-{\bf z}_{\star} \right\| _{2}^{2} \right)&      if\,\, 1\leqslant t_k\leqslant m,\\
			\mathbb{E} \left( \left\| {\bf x}_{k+1}-{\bf x}_{\star} \right\| _{2}^{2}+\left\| {\bf z}_{k+1}-{\bf z}_{\star} \right\| _{2}^{2} \right)\lesssim \tilde{\theta}_{k}\mathbb{E} \left( \left\| {\bf x}_k-{\bf x}_{\star} \right\| _{2}^{2}+\left\| {\bf z}_k-{\bf z}_{\star} \right\| _{2}^{2} \right)&      if\,\, m+1\leqslant t_k\leqslant m+n,\\
		\end{cases}
	\end{equation*}
where ``$\lesssim$" means the high order term $o \left(  \alpha \mathbb{E} \left\| {\bf x}_k-{\bf x}_{\star} \right\| _{2}^{2}+\tilde{\delta} \mathbb{E} \left\| {\bf z}_k-{\bf z}_{\star} \right\| _{2}^{2} \right)$ is omitted, and  $\tilde{\theta}_{k}<\max \big\{ \alpha ,\tilde{\delta} \big\} <1$.
\end{theorem}

\begin{remark}
As we only use a small portion of the rows of the data matrix $A$, we often have that $\hat{\beta}>\tilde{\beta}$ and $\tilde{\theta}_{k}>\theta _k$. That is, Algorithm \ref{alg6} may use more iterations than Algorithm \ref{alg5}. However, this does not mean that Algorithm \ref{alg6} definitely runs slower than Algorithm \ref{alg5}. Indeed, Algorithm \ref{alg6} is much cheaper than Algorithm \ref{alg5} per iteration, and it may use less CPU time in practice. One refers to Section 5 for numerical comparisons of the algorithms.
\end{remark}

\section{A practical stopping criterion for Kaczmarz-type methods}\label{Sec4}
\setcounter{equation}{0}
To the best of our knowledge, there are no practical stopping criteria for the randomized Kaczmarz-type methods till now. To fill-in this gap, we introduce a practical stopping criterion for the randomized Kaczmarz-type methods in this section. Let $\left\{ {\bf x}_i \right\} _{i=0}^{\infty}$ be an iterative sequence generated by {\it any} Kaczmarz-type methods, and let ${\bf x}_{\star}=A^{\dagger}{\bf b}$ be the least-squares solution of the linear system $A{\bf x}={\bf b}$. One scheme is to use the relative solution error (RSE) as the stopping criterion \cite{ref3, ref13, ref29, ref4}
\begin{equation}\label{3.1}
	RSE=\frac{\left\| {\bf x}_k-{\bf x}_{\star} \right\| _2}{\left\| {\bf x}_{\star} \right\| _2}.
\end{equation}
In \cite{ref7}, the absolute solution error (ASE) is used
\begin{equation}\label{3.2}
	ASE=\left\| {\bf x}_k-{\bf x}_{\star} \right\| _2.
\end{equation}
Obviously, both \eqref{3.1} and \eqref{3.2} are impractical since $x_{\star}$ is unknown in advance.
In \cite{ref8}, the relative residual (RRes) is utilized
\begin{equation*}
	RRes=\frac{\left\| {\bf b}-A{\bf x}_k \right\| _2}{\left\| {\bf b} \right\| _2}.
\end{equation*}
However, it is unsuitable for inconsistent linear system, and one has to calculate the residual ${\bf b}-A{\bf x}_k$. In other words, we have to access all the rows of $A$ in each iteration. In \cite{ref9}, Li and Wu suggest using the adjacent iterative solution error (AISE for short) as the stopping criterion
\begin{equation}\label{3.4}
	AISE=\frac{\left\| {\bf x}_k-{\bf x}_{k-1} \right\| _2}{\left\| {\bf b} \right\| _2}.
\end{equation}
We observe that it is free of $x_{\star}$ and there is no need to compute the residual vector in \eqref{3.4}. Unfortunately,
$\left\| {\bf x}_k-{\bf x}_{k-1} \right\| _2/\left\| {\bf b} \right\| _2$ converges to zero is only a {\it necessary} condition for $\left\{ {\bf x}_i \right\} _{i=0}^{\infty}$ converges to ${\bf x}_{\star}$. More precisely, 
\begin{equation*}
	\frac{\left| \left\| {\bf x}_k-{\bf x}_{\star} \right\| _2-\left\| {\bf x}_{k-1}-{\bf x}_{\star} \right\| _2 \right|}{\left\| {\bf b} \right\| _2}\leqslant \frac{\left\| {\bf x}_k-{\bf x}_{k-1} \right\| _2}{\left\| {\bf b} \right\| _2}\leqslant \frac{\left\| {\bf x}_k-{\bf x}_{\star} \right\| _2+\left\| {\bf x}_{k-1}-{\bf x}_{\star} \right\| _2}{\left\| {\bf b} \right\| _2}.
\end{equation*}


In this section, we will propose a practical stopping criterion for all the randomized Kaczmarz-type methods, and show its rationality theoretically. We need the following lemma.
\begin{Lemma}\label{Lem3.1}
Let $\left\{ {\bf x}_i \right\} _{i=0}^{\infty}$ be an iterative sequence obtained from a convergent Kaczmarz-type method. Then, there is a positive integer $L$ such that $\left\| {\bf x}_{\left( k+1 \right) L}-{\bf x}_{\star} \right\| _2<\left\| {\bf x}_{kL}-{\bf x}_{\star} \right\| _2, ~k=0,1,2,\ldots$
\end{Lemma}
\begin{proof}
We prove it by contradiction. Suppose that $\left\| {\bf x}_{\left( k+1 \right) L}-{\bf x}_{\star} \right\| _2\geqslant\left\| {\bf x}_{kL}-{\bf x}_{\star} \right\| _2, k=0,1,2,\ldots$, for all the positive integers $L$. If we take $L=1$, then
	\begin{equation*}
		\left\| {\bf x}_{k+1}-{\bf x}_{\star} \right\| _2\geqslant \left\| {\bf x}_k-{\bf x}_{\star} \right\| _2, \quad k=0,1,2,\ldots,
	\end{equation*}
     which contradicts to the fact $\left\{ {\bf x}_i \right\} _{i=0}^{\infty}$ is from a convergent Kaczmarz-type method.
\end{proof}

Let
\begin{equation*}
	g\left( \ell \right) =\left\| {\bf x}_{\ell}-{\bf x}_{\star} \right\| _2\,\,, \quad \ell=0,L,2L,\ldots
\end{equation*}
then the backward divided difference of $g\left( \ell \right)$ at $kL$ is defined as \cite{ref11}
    \begin{equation*}
    	\frac{\left| g\left( kL \right) -g\left( \left( k-1 \right) L \right) \right|}{kL-\left( k-1 \right) L}=\frac{\left| \left\| {\bf x}_{kL}-{\bf x}_{\star} \right\| _2-\left\| {\bf x}_{\left( k-1 \right) L}-{\bf x}_{\star} \right\| _2 \right|}{L}.
    \end{equation*}
From Lemma \ref{Lem3.1}, we see that $g\left( \ell \right)$ is a strictly monotonically decreasing discrete function whose infimum is 0, and
     \begin{equation*}
     	\lim_{k\rightarrow \infty} g\left( kL \right) =0.
     \end{equation*}
Thus, for any $\varepsilon >0$, there exists a positive integer $k_1>0$, such that $g\left( kL \right) \leqslant \frac{\varepsilon}{2},~\forall \,\,k\geqslant k_1$.
Next, we will show that there are $k_2>0$ and $\varepsilon _1<\varepsilon$, if $\left\| {\bf x}_{k_2L}-{\bf x}_{\left( k_2-1 \right) L} \right\| _2/L\leqslant \varepsilon _1$, then $g\left( k_2L \right) <\varepsilon$. Indeed, let $\varepsilon _1=\frac{\varepsilon}{2Lk_1}$, there exists $k_2>0$, such that
    \begin{align}
	    \frac{\left| g\left( k_2L \right) -g\left( \left( k_2-1 \right)L \right) \right|}{L}&=\frac{\left| \left\| {\bf x}_{k_2L}-{\bf x}_{\star} \right\| _2-\left\| {\bf x}_{\left( k_2-1 \right) L}-{\bf x}_{\star} \right\| _2 \right|}{L}\nonumber\\
	    &\leqslant \frac{\left\| {\bf x}_{k_2L}-{\bf x}_{\left( k_2-1 \right) L} \right\| _2}{L} \leqslant \frac{\varepsilon}{2Lk_1}.\label{3.8}
	\end{align}

On one hand, if $k_1\leqslant k_2$, then $g\left( k_2L \right) \leqslant \frac{\varepsilon}{2}<\varepsilon$. On the other hand, if $k_1>k_2$, then there is a scalar
    \begin{equation}\label{3.9}
    	a=\frac{g\left( k_1L \right) -g\left( k_2L \right)}{\left( k_1-k_2 \right) L},
    \end{equation}
such that
    \begin{equation}\label{3.10}
    	\frac{g\left( k_2L \right) -g\left( \left( k_2-1 \right) L \right)}{L}\leqslant a\leqslant \frac{g\left( k_1L \right) -g\left( \left( k_1-1 \right) L \right)}{L}< 0.
    \end{equation}
Combining the inequalities \eqref{3.8} and \eqref{3.10}, we arrive at
    \begin{equation*}
    	\left| a \right|\leqslant \frac{\left| g\left( k_2L \right) -g\left( \left( k_2-1 \right) L \right) \right|}{L}\leqslant \frac{\varepsilon}{2k_1L}.
    \end{equation*}
From \eqref{3.9}, we have $g\left( k_2L \right) =g\left( k_1L \right) -a\left( k_1-k_2 \right) L$. Hence,
    \begin{align*}
    	g\left( k_2L \right) &=\left| g\left( k_2L \right) \right| =\left|g\left( k_1L \right) -a\left( k_1-k_2 \right) L \right|\\
    	&\leqslant \left| g\left( k_1L \right) \right|+\left| a \right|\cdot \left| k_1-k_2 \right|L \leqslant \frac{\varepsilon}{2}+\frac{\varepsilon}{2k_1L}\left| k_1-k_2 \right|L\\
    	&=\frac{\varepsilon}{2}+\frac{k_1-k_2}{k_1}\cdot\frac{\varepsilon}{2} \leqslant \frac{\varepsilon}{2}+\frac{\varepsilon}{2}=\varepsilon.
    \end{align*}
In conclusion, we have the following theorem.
\begin{theorem}\label{Thm3.1}
Let $g\left( \ell \right) =\left\| {\bf x}_{\ell}-{\bf x}_{\star} \right\| _2~(\ell=0,L,2L,\ldots)$ be a strictly monotonically decreasing discrete function whose infimum is 0. For any $\varepsilon >0$, there are $k_2>0$ and $\varepsilon _1<\varepsilon$, if \[\frac{\left\| {\bf x}_{k_2L}-{\bf x}_{\left( k_2-1 \right) L} \right\|}{L}\leqslant \varepsilon _1,\] then $g\left( k_2 \right) <\varepsilon$.
\end{theorem}
\begin{remark}
Based on Theorem \ref{Thm3.1}, we propose the following stopping criterion for Kaczmarz-type methods
\begin{equation}\label{3.6}
	LISE=\frac{\left\| {\bf x}_{kL}-{\bf x}_{\left( k-1 \right) L} \right\| _2}{L}< tol,
\end{equation}
where $L>0$ is a user-provided number and $\it tol$ is the convergence tolerance.
\end{remark}

It is seen that the parameter $L$ is crucial to our stopping criterion. Next, we briefly discuss how to choose it in practice. Assume that
    \begin{equation*}
    	\mathbb{E} \left\| {\bf x}_{k+1}-{\bf x}_{\star} \right\| _{2}^{2}\leqslant c\mathbb{E} \left\| {\bf x}_k-{\bf x}_{\star} \right\| _{2}^{2},
    \end{equation*}
where $0<c<1$, and
    \begin{equation}\label{3.12}
    	\mathbb{E} \left\| {\bf x}_{\left( k+1 \right) L}-{\bf x}_{\star} \right\| _{2}^{2}\leqslant c^L\mathbb{E} \left\| {\bf x}_{kL}-{\bf x}_{\star} \right\| _{2}^{2}.
    \end{equation}
Therefore,
    \begin{align*}
    	\mathbb{E} \left\| {\bf x}_{kL}-{\bf x}_{\star} \right\| _{2}^{2}&=\mathbb{E} \left\| {\bf x}_{kL}-{\bf x}_{\left( k+1 \right) L}+{\bf x}_{\left( k+1 \right) L}-{\bf x}_{\star} \right\| _{2}^{2} \leqslant \mathbb{E} \left( \left\| {\bf x}_{\left( k+1 \right) L}-{\bf x}_{kL} \right\| _2+\left\| {\bf x}_{\left( k+1 \right) L}-{\bf x}_{\star} \right\| _2 \right) ^2\\
    	&=\mathbb{E} \left( \left\| {\bf x}_{\left( k+1 \right) L}-{\bf x}_{kL} \right\| _{2}^{2}+2\left\| {\bf x}_{\left( k+1 \right) L}-{\bf x}_{kL} \right\| _2\left\| {\bf x}_{\left( k+1 \right) L}-{\bf x}_{\star} \right\| _2+\left\| {\bf x}_{\left( k+1 \right) L}-{\bf x}_{\star} \right\| _{2}^{2} \right)\\
    	&=\mathbb{E} \left\| {\bf x}_{\left( k+1 \right) L}-{\bf x}_{kL} \right\| _{2}^{2}+2\mathbb{E} \left( \left\| {\bf x}_{\left( k+1 \right) L}-{\bf x}_{kL} \right\| _2\left\| {\bf x}_{\left( k+1 \right) L}-{\bf x}_{\star} \right\| _2 \right) +\mathbb{E} \left\| {\bf x}_{\left( k+1 \right) L}-{\bf x}_{\star} \right\| _{2}^{2}\\
    	&\leqslant \mathbb{E} \left\| {\bf x}_{\left( k+1 \right) L}-{\bf x}_{kL} \right\| _{2}^{2}+2\mathbb{E} \left( \left\| {\bf x}_{\left( k+1 \right) L}-{\bf x}_{kL} \right\| _2\left\| {\bf x}_{\left( k+1 \right) L}-{\bf x}_{\star} \right\| _2 \right) +c^L \mathbb{E} \left\| {\bf x}_{kL}-{\bf x}_{\star} \right\| _{2}^{2}.
    \end{align*}

As a result,
    \begin{equation*}
    	\left( 1-c^L \right) \mathbb{E} \left\| {\bf x}_{kL}-{\bf x}_{\star} \right\| _{2}^{2}\leqslant \mathbb{E} \left\| {\bf x}_{\left( k+1 \right) L}-{\bf x}_{kL} \right\| _{2}^{2}+2\mathbb{E} \left( \left\| {\bf x}_{\left( k+1 \right) L}-{\bf x}_{kL} \right\| _2\left\| {\bf x}_{\left( k+1 \right) L}-{\bf x}_{\star} \right\| _2 \right),
    \end{equation*}
and
    \begin{equation}\label{3.13}
    	\mathbb{E} \left\| {\bf x}_{kL}-{\bf x}_{\star} \right\| _{2}^{2}\leqslant \frac{1}{1-c^L}\mathbb{E} \left\| {\bf x}_{\left( k+1 \right) L}-{\bf x}_{kL} \right\| _{2}^{2}+\frac{2}{1-c^L}\mathbb{E} \left( \left\| {\bf x}_{\left( k+1 \right) L}-{\bf x}_{kL} \right\| _2\left\| {\bf x}_{\left( k+1 \right) L}-{\bf x}_{\star} \right\| _2 \right).
    \end{equation}
From \eqref{3.12} and \eqref{3.13}, it follows that
    \begin{equation*}
    	\mathbb{E} \left\| {\bf x}_{\left( k+1 \right) L}-{\bf x}_{\star} \right\| _{2}^{2}\leqslant \frac{c^L}{1-c^L}\mathbb{E} \left\| {\bf x}_{\left( k+1 \right) L}-{\bf x}_{kL} \right\| _{2}^{2}+\frac{2c^L}{1-c^L}\mathbb{E} \left( \left\| {\bf x}_{\left( k+1 \right) L}-{\bf x}_{kL} \right\| _2\left\| {\bf x}_{\left( k+1 \right) L}-{\bf x}_{\star} \right\| _2 \right).
    \end{equation*}
This implies that we need to use a sufficiently large $L$ in practical calculations. One can use, say, $L=50$ or $80$ in practice.



\section{Numerical Experiments}\label{sec5}

In this section, we perform numerical experiments to show the efficiency of our algorithms as well as the effectiveness of the proposed stopping criterion. All the numerical experiments are run on a Hp workstation with 20 cores double Intel(R)Xeon(R) E5-2640 v3 processors, with CPU 2.60 GHz and RAM 256 GB. The operation system is 64-bit Windows 10. All the numerical results are obtained from using MATLAB 2018b. In order to show the efficiency of our proposed algorithms, we compare Algorithm \ref{alg5} and Algorithm \ref{alg6} with the following state-of-the-art Kaczmarz-type methods for solving  large-scale (dense or sparse) inconsistent linear systems:\\
$\bullet$ {\bf REK \cite{ref5}}: The randomized extended Kaczmarz method.\\
$\bullet$ {\bf PREK \cite{ref13}}: The partially randomized extended Kaczmarz method.\\
$\bullet$ {\bf GRAK \cite{ref6}}: The greedy randomized augmented Kaczmarz method.\\
$\bullet$ {\bf RGRAK \cite{ref12}}: The relaxed greedy randomized augmented Kaczmarz method, in which the relaxation parameter is set to be 1 \cite{ref12}.



In the tables below, we denote by ``IT" the number of iterations, by ``CPU" the CPU time in seconds, and by ``RSE" the relative error defined as
\begin{equation}
{\rm RSE}=\frac{\left\| {\bf x}_k-{\bf x}_{\star} \right\| _2}{\left\| {\bf x}_{\star} \right\| _2},
\end{equation}
where ${\bf x}_{\star}=A^{\dagger}{\bf b}$ is the least-squares solution in 2-norm of the inconsistent linear systems \eqref{1.1}, which is generated by the MATLAB function $pinv.m$. All the numerical results are the means from 10 runs.


\subsection{Effectiveness of the stopping criterion \eqref{3.6}}\label{Sec5.1}

In this section, we show effectiveness of our stopping criterion for large inconsistent linear systems. In this subsection, the test matrix $A$ is synthetic and is generated by using the MATLAB built-in function $randn(m,n)$. Similar to \cite{ref6}, the right-hand side is chosen as ${\bf b}=A{\bf x}_{\ast}+{\bf r}$, where ${\bf x}_{\ast}\in \mathbb{R} ^n$ is one of the least-squares solution of the inconsistent linear systems \eqref{1.1}, and ${\bf r}\in \mathbb{R} ^m$ is a nonzero vector in the null space of the matrix $A^H$ generate by the MATLAB function $null.m$.
Two examples are performed in this subsection.
%

$\bullet$ In this first example, we run the six algorithms using their own stopping criteria. More precisely, in the REK method, we use \cite{ref5}
\begin{equation}\label{5.1}
	\frac{\left\| A{\bf x}_k-\left( {\bf b}-{\bf z}_k \right) \right\| _2}{\left\| A \right\| _F\left\| {\bf x}_k \right\| _2}\leqslant tol\quad {\rm and} \quad\frac{\left\| A^H{\bf z}_k \right\| _2}{\left\| A \right\| _{F}^{2}\left\| {\bf x}_k \right\| _2}\leqslant tol
\end{equation}
as the stopping criterion; in the PREK method, we use \cite{ref13}
\begin{equation}\label{5.2}
    \frac{\left\| {\bf x}_k-{\bf x}_{\star} \right\| _{2}^{2}}{\left\| {\bf x}_{\star} \right\| _{2}^{2}}\leqslant tol
\end{equation}
as the stopping criterion; in the GRAK method and the RGRAK method, we use \cite{ref6, ref12}
\begin{equation}\label{5.3}
	\frac{\left\| {\bf x}_k-{\bf x}_{\star} \right\| _{2}^{2}+\left\| {\bf z}_k-{\bf z}_{\star} \right\| _{2}^{2}}{\left\| {\bf x}_{\star} \right\| _{2}^{2}+\left\| {\bf b}_{\mathcal{R} \left( A \right)} \right\| _{2}^{2}}\leqslant tol
\end{equation}
as the stopping criteria, where ${\bf z}_{\star}={\bf b}_{\mathcal{R} \left( A \right) ^{\bot}}$ is computed by ${\bf b}_{\mathcal{R} \left( A \right) ^{\bot}}={\bf b}-A{\bf x}_{\star}$. In Algorithm \ref{alg5} and Algorithm \ref{alg6}, we make use of
\begin{equation}\label{55}
	LISE =\frac{\left\| \tilde{{\bf x}}_{kL}-\tilde{{\bf x}}_{\left( k-1 \right) L} \right\| _2}{L}\leqslant tol
\end{equation}
as the stopping criteria, where $L=400$, $\tilde{{\bf x}}_{kL}=\left[ {\bf z}_{kL}^H, {\bf x}_{kL}^H \right]^H$, refer to \eqref{eq22}, and the convergence tolerance is chosen as $tol=10^{-4}$. Table \ref{T2} lists the numerical results, where the values in bold are the best one.


\begin{table}[H]
	\centering
	\caption{{\bf Section 5.1.} Example 1: Numerical results of the algorithms using their own stopping criteria \eqref{5.1}--\eqref{55}; $A=randn(m,n)$, and $tol=10^{-4}$. The sampling ratio is chosen as $\eta=0.01$ in Algorithm \ref{alg6}.}
	\resizebox{.99\columnwidth}{!}{
	{\footnotesize\begin{tabular}{|c|c|c|c|c|c|c|}
			\hline
			\multicolumn{2}{|c|}{$m\times n$}   & {$5000\times1000$} & {$5000\times1500$} & {$5000\times2000$} & {$5000\times2500$} & {$5000\times3000$} \\ \hline
			\multirow{2}{*}{REK}   & IT  & 19916     & 37204.9     & 64410.7     & 109528.8     & 185659.2     \\ \cline{2-7}
			& CPU & 161.9   & 375.7     & 646.4   & 1362.3   & 2618.0    \\
			\cline{2-7}
			& RSE & 4.38e-3   & 6.53e-3     & 9.44e-3   & 1.34e-2   & 1.98e-2    \\ \hline
            \multirow{2}{*}{PREK}   & IT  & 13191.5     & 24724.2     & 47252.7     & 87730.9     & 169481.5     \\ \cline{2-7}
			& CPU & {\bf 3.5}   & {\bf 7.3}     & {\bf 14.9}   & {\bf 31.4}   & {\bf 64.3}    \\
			\cline{2-7}
			& RSE & 1.00e-2   & 1.00e-2     & 1.00e-2   & 1.00e-2   & 1.00e-2    \\ \hline
			\multirow{2}{*}{GRAK}  & IT  & 3003.5     & 5645.6      & { \bf 10181.9}      & 16879.9      & 29117.3     \\ \cline{2-7}
			& CPU & 17.1     & 49.9    & 84.7    & 139.9    & 282.1   \\
			\cline{2-7}
			& RSE & 4.79e-1   & 4.50e-1     & 4.16e-1   & 3.94e-1   & 3.61e-1    \\ \hline
            \multirow{2}{*}{RGRAK}  & IT  & {\bf 2977.3}     & {\bf 5601.4}      & 10206      & {\bf 16702.6}      & {\bf 28927.7}     \\ \cline{2-7}
			& CPU & 15.7     & 47.6    & 84.1    & 139.6    & 283.5   \\
			\cline{2-7}
			& RSE & 4.90e-1   & 4.73e-1     & 4.55e-1   & 4.44e-1   & 4.10e-1    \\ \hline
			\multirow{2}{*}{\bf Algorithm \ref{alg5}}  & IT  & 9600     & 19200      & 36880      & 65480      & 127200     \\ \cline{2-7}
			& CPU & 61.1     & 117.9    & 253.6    & 581.5    & 1356.9   \\
			\cline{2-7}
			& RSE & 7.52e-4   & 1.17e-3     & 1.73e-3   & 2.22e-3   & 3.19e-3    \\ \hline
			\multirow{2}{*}{\bf Algorithm \ref{alg6}}  & IT  & 10240     & 19920     & 37960     & 68080      & 131360      \\ \cline{2-7}
			& CPU & 14.0     & 38.9    & 107.5    & 250.7   & 748.6   \\
			\cline{2-7}
			& RSE & {\bf 6.61e-4}   & {\bf 1.06e-3}     & {\bf 1.63e-3}   & {\bf 2.09e-3}   & {\bf 2.99e-3}    \\ \hline
			
	\end{tabular}}\label{T2}}
\end{table}

It is seen from Table \ref{T2} that the PREK method and the RGRAK method may require less CPU time and fewer number of iterations by using their own stopping criterion. However, the computed solution obtained from the two methods are unreliable, and the ``real" relative errors of the approximations from the two methods can be large. For instance, the RSE values of the GRAM and RGRAK methods are only in the order of $\mathcal{O}(10^{-1})$, the RSE values of the REK and PREK methods are only in the order of $\mathcal{O}(10^{-2})$ or $\mathcal{O}(10^{-3})$, while those of Algorithm \ref{alg5} and Algorithm \ref{alg6} are in the order of $\mathcal{O}(10^{-3})$ or $\mathcal{O}(10^{-4})$.

Indeed, we find that the denominator of the stopping criterion \eqref{5.3} can be quite large, which leads to an early termination of the GRAM and RGRAK methods even when ${\bf x}_k$ is far from being a good approximation to ${\bf x}_{\star}$. For example, for the matrix $A$ with $m=5000$ and $n=1000$, the values of $\left\| {\bf x}_{\star} \right\| _{2}^{2}+\left\| {\bf b}_{\mathcal{R} \left( A \right)} \right\| _{2}^{2}$ can be over $5\times 10^6$.
Moreover, both \eqref{5.2} and \eqref{5.3} are impractical as $x_{\star}$ is unknown in advance, and we have to perform two matrix-vector product with respect to $A$ and $A^H$ for computing \eqref{5.1}. That is to say, one has to access all the rows and columns of $A$. Thus, compared with some available stopping criteria, our stopping criterion \eqref{55} is effective and is more appropriate as a stopping criterion for Kaczmarz-type methods.

\begin{table}[H]
	\centering
	\caption{{\bf Section 5.1.} Example 2: Numerical results of the algorithms using the proposed stopping criterion \eqref{55} with $L=400$; $A=randn(m,n)$ and $tol=10^{-4}$. The sampling ratio is chosen as $\eta=0.01$ in Algorithm \ref{alg6}.}
	\resizebox{.99\columnwidth}{!}{
	{\footnotesize\begin{tabular}{|c|c|c|c|c|c|c|}
			\hline
			\multicolumn{2}{|c|}{$m\times n$}   & {$5000\times1000$} & {$5000\times1500$} & {$5000\times2000$} & {$5000\times2500$} & {$5000\times3000$} \\ \hline
			\multirow{2}{*}{REK}   & IT  & 27680     & 54360     & 98480     & 173520     & 306360     \\ \cline{2-7}
			& CPU & 9.5   & 21.3     & 41.4   &73.3   & 139.1    \\
			\cline{2-7}
			& RSE & {\bf 4.43e-4}   & {\bf 7.04e-4}     & {\bf 9.94e-4}   & {\bf 1.46e-3}   & {\bf 2.28e-3}    \\ \hline
            \multirow{2}{*}{PREK}   & IT  & 17600     & 33360     & 60520     & 106200     & 186120     \\ \cline{2-7}
			& CPU & {\bf 4.7}   & {\bf 9.9}     & {\bf 19.4}   &{\bf 36.7}   & {\bf 66.5}    \\
			\cline{2-7}
			& RSE & 2.17e-3   & 2.87e-3     & 3.49e-3   & 4.52e-3   & 5.84e-3    \\ \hline
			\multirow{2}{*}{GRAK}  & IT  & 11200     & 23040      & 43760      & 81960      & 146280     \\ \cline{2-7}
			& CPU & 73.9     & 230.9    & 455.4    & 1013.4    & 1889.4   \\
			\cline{2-7}
			& RSE & 1.01e-3   & 1.63e-3     & 2.22e-3   & 3.04e-3   & 4.13e-3    \\ \hline
            \multirow{2}{*}{RGRAK}  & IT  & 11040     & 22640      & 42520      & 78640      & 140800     \\ \cline{2-7}
			& CPU & 73.5     & 184.1    & 501.6    & 1123.3    & 1961.4   \\
			\cline{2-7}
			& RSE & 1.02e-3   & 1.76e-3     & 2.59e-3   & 3.61e-3   & 5.58e-3    \\ \hline
			\multirow{2}{*}{\bf Algorithm \ref{alg5}}  & IT  & {\bf 9600}     & {\bf 19200}      & {\bf 36040}   & {\bf 66440}      & {\bf 116600}     \\ \cline{2-7}
			& CPU & 57.1     & 96.1    & 289.2    & 592.9    & 1119.5   \\
			\cline{2-7}
			& RSE & 7.77e-4   & 1.19e-3     & 1.66e-3   & 2.23e-3   & 2.99e-3    \\ \hline
			\multirow{2}{*}{\bf Algorithm \ref{alg6} }  & IT  & 10120     & 20320     & 37720     & 68920      & 120760      \\ \cline{2-7}
			& CPU & 14.0     & 40.6    & 110.2    & 223.3   & 713.5   \\
			\cline{2-7}
			& RSE & 6.99e-4   & 1.10e-3     & 1.54e-3   & 2.09e-3   & 2.88e-3    \\ \hline
			
	\end{tabular}}\label{T3}}
\end{table}


$\bullet$ In view of the numerical results given in the first example, in the second example, we run the six algorithms using the proposed stopping criterion \eqref{3.6} with $L=400$ and $tol=10^{-4}$.  We generate another five test matrices by using the MATLAB built-in function $randn(m,n)$.
The numerical results are given in Table \ref{T3}, where the best results are in bold.

Some remarks are in order. First, comparing the $RSE$ values of the GRAK, RGRAK, REK and PREK methods in Tabel \ref{T2} and Tabel \ref{T3}, we find that the approximate solutions got from these four method are more accurate than before when using \eqref{5.1} as the stopping criterion. This means that the proposed stopping criterion is feasible and it is more effective than the available ones.
Second, we observe from Table \ref{T3} and Table \ref{T2} that REK runs much faster when using \eqref{55} as the stopping criterion. This is due to the fact that one has to perform two matrix-vector product with respect to $A$ and $A^H$ when evaluating \eqref{5.1}. In fact, the computational cost of matrix-vector products is high when $A$ is large and dense. Third, it is seen from the table that for this synthetic problem, REK enjoys the highest accuracy and PREK runs the fastest, while Algorithm \ref{alg5} needs the fewest iterations. This implies that the convergence speed of Algorithm \ref{alg5} is the highest among the six. Of cause, an algorithm with least iteration numbers does not mean that it runs the fastest. Indeed, the framework of REK and PREK is different from those of the other four algorithms. Fourth, we see that Algorithm \ref{alg5} and Algorithm \ref{alg6} outperform GRAK and RGRAK significantly, both in terms of CPU time and iteration numbers. Specifically, Algorithm \ref{alg5} needs the fewest iterations and Algorithm \ref{alg6} uses the least CPU time. This is because Algorithm \ref{alg6} only utilize a small portion of rows of the matrix $\tilde{A}$ in each iteration.

\subsection{Efficiency of the proposed algorithms}\label{Sec5.2}

In this section, we perform numerical experiments to show the superiority of the proposed algorithms over many state-of-the-art Kaczmarz-type methods for large and inconsistent linear systems. Based on the numerical experiments made in Section \ref{Sec5.1} and for the sake of justification, we make use of \eqref{3.6} as the stopping criterion for all the algorithms, with $L=400$. In order to demonstrate the superiority of Algorithm \ref{alg6} over the GRAK method, we consider the value of
\begin{equation}\label{5.6}
   {\rm Speed \text{-} up}=\frac{\rm CPU~time~of~GRAK}{\rm CPU~time~of~Algorithm~\ref{alg6}}.
\end{equation}
There are three numerical experiments in this section.

\begin{table}[H]
	\caption{{\bf Section 5.2.} Example 1: Test matrices for solving large inconsistent linear systems.}
	\centering
	{\footnotesize\begin{tabular}{|c|c|c|c|}
			\hline
			Matrix ($A$) & Size ($m\times n$) & Nnz & Background \\\hline
			$abtaha1\footnotemark[2]$ & {$14596\times209$}& 51307 & Combinatorial Problem \\
			$abtaha2\footnotemark[2]$ & {$37932\times331$} & 137228 & Combinatorial Problem \\
			$sls\footnotemark[2]$ & {1748122$\times62729$} & 6804304 & Least Squares Problem \\
			$California\footnotemark[2]$ & {$9664\times9664$} & 16150 & Directed Graph \\
			$stat96v5\footnotemark[2]$ & {$2307\times75779$} & 233921 & Linear Programming Problem \\
            $GL7d25\footnotemark[2]$ & {$2798\times21074$} & 81671 & Combinatorial Problem \\
			\hline
	\end{tabular}}\label{T1}
\end{table}

$\bullet$ In the first experiment, we show numerical behavior of the proposed algorithms on some large sparse inconsistent linear systems.
The test matrices are from the University of Florida Sparse Matrix Collection \footnote[2]{https://sparse.tamu.edu/}. The details of these data matrices are listed in Table \ref{T1}. In this experiment, the construction of the right-hand side $\bf b$ corresponding to the matrices $abtaha1$, $abtaha2$, $California$, $stat96v5$ and $GL7d25$ is the same way as in Section \ref{Sec5.1}, while for the matrix $sls$, it is generated by using the MATLAB built-in function $randn$. The convergence tolerance is chosen as $tol=10^{-4}$, and the sampling ratio in Algorithm \ref{alg6} is set to be $\eta =0.01$. Table \ref{T4} lists the numerical results, where the values in bold are the best one.

\begin{table}[H]
	\centering
	\caption{{\bf Section \ref{Sec5.2}.} Example 1: Numerical results of the six algorithms on some sparse matrices from the University of Florida Sparse Matrix Collection, $tol=10^{-4}$ and $L=400$. The sampling ratio is chosen as $\eta=0.01$ in Algorithm \ref{alg6}. Here ``/" means that the number of iterations exceeds 1000000 or the CPU time exceeds 24 hours.}
	\resizebox{.99\columnwidth}{!}{
	{\footnotesize\begin{tabular}{|c|c|c|c|c|c|c|c|}
			\hline
			\multicolumn{2}{|c|}{Matrix\& Size}   & \makecell[c]{$abtaha1$ \\{$14596\times209$}} & \makecell[c]{$abtaha2$ \\{$37932\times331$}} & \makecell[c]{$sls$\\{1748122$\times62729$}} & \makecell[c]{$California$\\$9664\times9664$}& \makecell[c]{$stat96v5$\\$2307\times75779$} & \makecell[c]{$GL7d25$\\$2798\times21074$}\\ \hline
			\multirow{2}{*}{REK}   & IT  & 49640     & 64760     & /     & 476680     & 69440   & 103840  \\ \cline{2-8}
			& CPU & 21.8   & 51.6     & /   & 396.7   & 273.4   &152.8 \\ \hline
			\multirow{2}{*}{PREK}   & IT  & 36640     & 41800     & /     & /     & 100880   &/  \\ \cline{2-8}
			& CPU & 14.1   & 29.7     & /   & /   & 175.6   &/ \\ \hline
			\multirow{2}{*}{GRAK}  & IT  & 38880     & 56200      & /      & 300320      & 18680   &73880  \\ \cline{2-8}
			& CPU & 35.1     & 99.4    & /    & 418.2    & 119.6  &190.9 \\ \hline
			\multirow{2}{*}{RGRAK}   & IT  & 33680     & 47440     & /     & 398080    & 18720   & 93240 \\ \cline{2-8}
			& CPU & 29.5   & 85.5     & /   & 531.9   & 118.2   & 202.9 \\ \hline
			\multirow{2}{*}{\bf Algorithm \ref{alg5}}  & IT  & {\bf 25680}     & {\bf 32560}      & /      & {\bf 223280}      & {\bf 11800}    & {\bf 54440} \\ \cline{2-8}
			& CPU & {\bf 11.9}     & {\bf 29.2}    & /    & 165.7    & 33.4  & 70.1 \\ \hline
			\multirow{2}{*}{\bf Algorithm \ref{alg6} }  & IT  & 33880     & 45920     & 406640     & 251760      & 12680    & 59040  \\ \cline{2-8}
			& CPU & 15.6     & 45.1    & 59120.4    & {\bf 161.2}   & {\bf 31.8}  & {\bf 68.5} \\ \hline
            \multicolumn{2}{|c|}{\bf Speed-up}   &2.25   &2.20   &/   &2.59   &3.76   &2.79    \\ \hline
	\end{tabular}}\label{T4}}
\end{table}

Some remarks are in order. First, it is seen from Table \ref{T4} that Algorithm \ref{alg5} and Algorithm \ref{alg6} are far superior to the other algorithms in terms of CPU time and number of iterations. In terms of the values of ``Speed-up", we observe from Table \ref{T4} that Algorithm \ref{alg6} is about 2 to 3 times faster than the GRAK method, and the improvements are impressive. Specifically, for the large matrix $sls$, all the algorithms fail to converge within 24 hours, except for Algorithm \ref{alg6}. Therefore, we benefit from the strategy of simple sampling and there is no need to compute the whole residual vector during iterations. Thus, Algorithm \ref{alg6} is competitive for solving large and sparse inconsistent linear systems. Notice that the PREK method does not work for the two matrices $California$ and $GL7d25$, either. Indeed, it is required that the coefficient matrix $A$ is tall and
of full column rank in the PREK method \cite[Theorem 3.1]{ref13}, while the two matrices are not of full column rank.

$\bullet$ In the second experiment, we consider the problem of tomographic image reconstruction and run the six algorithms on a 2-D parallel-beam tomography problem. In this experiment, a sparse matrix $A\in \mathbb{R} ^{m\times n}$ and an ``exact solution" ${\bf x}_{\star}\in \mathbb{R} ^n$ is generated by using the function $paralleltomo(N, \theta, p)$ in the AIR tool box \cite{ref14}, with $N=60$, $\theta =0:1:178\degree$ and $p=125$. The matrix $A$ is of size $22375\times 3600$, and we set the right-hand side ${\bf b}=A{\bf x}_{\star}+{\bf r}$, where the noise ${\bf r}\in \mathbb{R} ^m$ is a nonzero vector in the null space of the matrix $A^H$ generate by the MATLAB function $null$.

\begin{figure}[H]\label{Fig4.1}
	\begin{minipage}{0.99\linewidth}
		\centerline{\includegraphics[width=15cm]{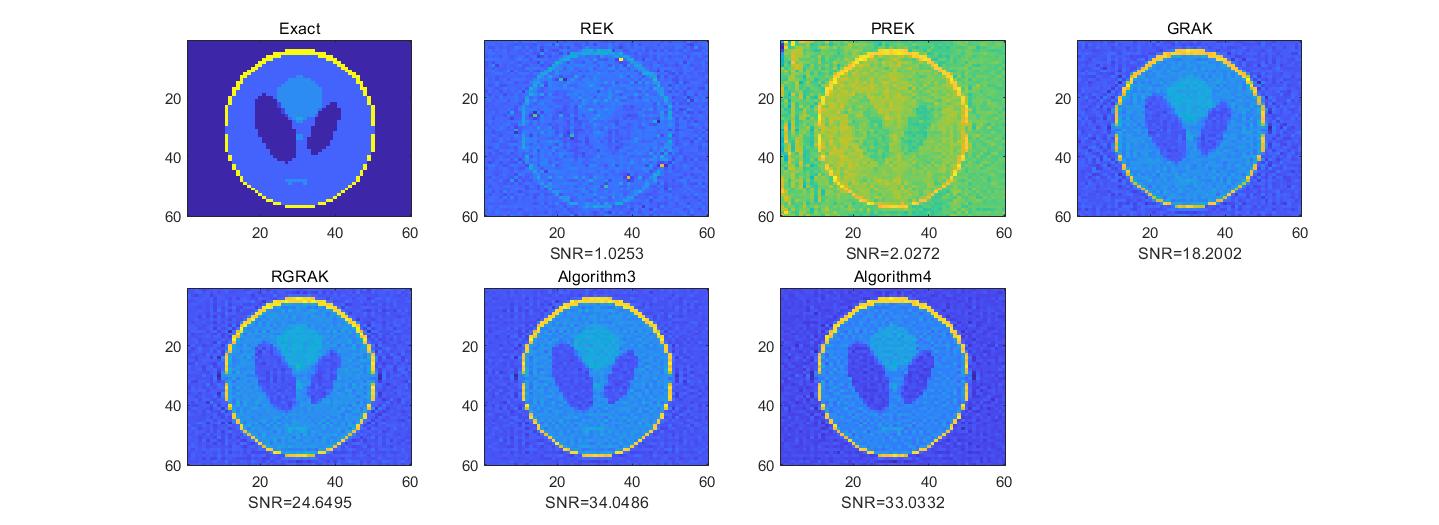}}
	\end{minipage}
	\caption{Example 2: The original image (Exact), the reconstructed ones, and the values of SNR
obtained from the six algorithms REK, PREK, GRAK, RGRAK, Algorithm \ref{alg5} and Algorithm \ref{alg6} $\left(\eta =0.01\right)$. All the algorithms were run for 20000 iterations.}
	\label{fig:res}
\end{figure}

In this experiment, we solve the inconsistent linear system $A{\bf x}={\bf b}$, and all the algorithms are run for 20000 iterations.
Consider the signal-noise ratio (SNR) \cite{ref41}
\begin{equation}\label{5.4}
  SNR= \frac{\sum_{i=1}^n{{\bf x}_{i}^{2}}}{\sum_{i=1}^n{\left( {\bf x}_i-\hat{{\bf x}}_i \right) ^2}},
\end{equation}
where ${\bf x}$ is the original clean signal (i.e., the ``exact solution" ${\bf x}_{\star}$), and $\hat{{\bf x}}$ is the denoised signal (i.e., the computed solution obtained from one of the six algorithms). Notice that the larger the $SNR$, the better the denoising effect of an algorithm. The original image (Exact), the reconstructed ones as well as the values of SNR are shown in Figure \ref{fig:res}. It is obviously to see from the figures and the values of SNR that Algorithm \ref{alg5} and Algorithm \ref{alg6} can remove the noise and restore the real image efficiently, and our algorithms are superior to the other methods in terms of SNR.

$\bullet$ Spectral Regression Discriminant Analysis (SRDA) is one of the most popular methods for large-scale discriminant analysis \cite{SRDA}. By using the tool of spectral graph analysis, SRDA casts discriminant analysis into a regression framework which facilitates both efficient computation and the use of regularization techniques. Specifically, SRDA only needs to solve a set of least squares problems and there is no eigenvector computation involved, which can save both the running time and the memory compared with some classical linear discriminant analysis methods.

Let $X=\left[ {\bf x}_1,{\bf x}_2,\dots ,{\bf x}_n \right] \in \mathbb{R} ^{d\times n}$ be a set of $n$ training samples in a $d$-dimensional feature space. Assume that the data matrix is partitioned into $k$ classes as $X=\left[ X_1,\dots ,X_k \right]$, where $X_j$ is the $j$-th set with $n_j$ being the number of samples. Denote by ${\bf c}_j$ the centroid vector of $X_j$, and by ${\bf c}$ the global centroid vector of the training data. Denote by
\begin{equation*}
  Y=\left[ {\bf y}_1,{\bf y}_2,\dots ,{\bf y}_k \right] =\left[ \begin{matrix}
           {\bf 1}_{n_1}&               {\bf 0}&               \cdots&               {\bf 0}\\
           {\bf 0}&               {\bf 1}_{n_2}&               \cdots&               {\bf 0}\\
           \vdots&           \vdots&               \ddots&           \vdots\\
           {\bf 0}&               {\bf 0}&                     \cdots&           {\bf 1}_{n_k}\\
  \end{matrix} \right] \in \mathbb{R} ^{n\times k},
\end{equation*}
 where ${\bf 1}_{n_i}$ and ${\bf 1}_{n}$ are the one vector of size $n_i$ and $n$, respectively, and $I_n$ is the $n$-by-$n$ identity matrix. Let $\hat{Y}=\left( I_n-{\bf 1}_n{\bf 1}_{n}^{T} \right) Y$ and $\bar{Y}\in \mathbb{R} ^{n\times \left( k-1 \right)}$ be an orthonormal basis of $span\left\{ \hat{Y} \right\}$.

Let $\hat{X}=\left[ X^T,{\bf 1}_n \right] ^T$. In essence, the SRDA method resorts to the following least squares problems \cite{SRDA, ref50}:
\begin{equation}\label{5.5}
  \hat{a}_i=\mathop{\arg \min} \limits_{a_i\in \mathbb{R} ^d} \left\{ \left\| \hat{X}^Ta_i-{\bf y}_i \right\| _2\right\}, \quad i=1,2,\dots,k-1,
\end{equation}
where ${\bf y}_i$ is the $i$-th column vector of $\bar{Y},~i=1,2,\dots ,k-1$.

Without loss of generality, we solve the problem as $i=1$ by using Algorithm \ref{alg5} and Algorithm \ref{alg6}, where the sampling ratio is chosen as $\eta=0.001$ in the latter algorithm. We compare the two proposed methods with the REK method, the PREK method, the GRAK method, and the RGRAK method. In the stopping criterion \eqref{3.6}, we use $tol=10^{-6}$ and $L=400$. We run the algorithms on five dense databases, including the $MNIST$ dataset\footnote[3]{http://yann.lecun.com/exdb/mnist/}, the $UCI$ dataset\footnote[4]{http://archive.ics.uci.edu/dataset/54/isolet}, the $USPS\_$ dataset\footnote[5]{https://hyper.ai/datasets/16041}\cite{ref51}, the $CIFAR$-10 dataset\footnote[6]{http://www.cs.toronto.edu/~kriz/cifar.html}, as well as the YouTube Faces dataset\footnote[7]{http://www.cs.tau.ac.il/ wolf/ytfaces/}. In all the experiments, we randomly choose  $n=70\%\cdot N$ samples as the training set, where $N$ is the total number of samples. The details of the data sets are listed in Table \ref{den}, and Table \ref{T5} lists the numerical results obtained.
\begin{table}[H]
	\caption{{\bf Section 5.2.} Example 3: Details of the databases.}
	\centering
    \resizebox{.99\columnwidth}{!}{
	{\footnotesize\begin{tabular}{|c|c|c|c|c|}
			\hline
			Datasets & Dimensionality ($d$) & Number of total samples ($N$) & Background    &Type \\\hline
			$MNIST\footnotemark[3]$ &784 & 70000 & Handwritten digits recognition  & dense \\
			$ISOLET\footnotemark[4]$ & 617 & 7797 & Classification & dense \\
			$USPS\_\footnotemark[5]$ & 256 & 9598 & Handwritten digits recognition & dense \\
			$CIFAR-10\footnotemark[6]$ & 3072 & 50000 & Object recognition & dense \\
			$YouTubeFaceArrange\_64\times64\footnotemark[7]$ & 4096 & 370319 & Face recognition & dense \\
			\hline
	\end{tabular}}\label{den}}
\end{table}

\begin{table}[H]
	\centering
	\caption{{\bf Section 5.2.2.} Example 3: Numerical results of the algorithms for solving the least squares problem \eqref{5.5} as $i=1$. Here "/" means the number of iterations exceeds 1000000 or the CPU time exceeds 24 hours.}
	\resizebox{.99\columnwidth}{!}{
	{\footnotesize\begin{tabular}{|c|c|c|c|c|c|c|}
			\hline
            \multicolumn{2}{|c|}{Matrix\& Size}   & \makecell[c]{$MNIST$ \\{$70000\times784$}} & \makecell[c]{$ISOLET$ \\{$7797\times617$}} & \makecell[c]{$USPS\_$\\{9598$\times256$}} &
            \makecell[c]{$CIFAR-10$\\{50000$\times3072$}} &
            \makecell[c]{$YouTubeFaceArrange\_64\times64$\\$370319\times4096$}\\ \hline
			\multirow{2}{*}{REK}   & IT  & 166800   & 246440     & 367000     & /     & /     \\ \cline{2-7}
			& CPU & 237.6   & 83.3     &111.9   & /   & /    \\ \hline
            \multirow{2}{*}{PREK}   & IT  & /   & 374680     & 162600     & /     & /     \\ \cline{2-7}
			& CPU & /   & 101.6     &40.6   & /   & /    \\ \hline
			\multirow{2}{*}{GRAK}  & IT  & 19200     & 124840      & 188480      & 465080   & /     \\ \cline{2-7}
			& CPU & 419.9     & 312.5    & 153.1    & 30428.1   & /   \\ \hline
            \multirow{2}{*}{RGRAK}  & IT  & 17760     & 88600      & 98440      & 130800   & /     \\ \cline{2-7}
			& CPU & 350.2     & 217.9    & 79.9    & 14115.9   & /   \\ \hline
			\multirow{2}{*}{\bf Algorithm \ref{alg5}}  & IT  & {\bf 17200}     & {\bf 82200}      & {\bf 72320}      & {\bf 130760}   & /     \\ \cline{2-7}
			& CPU & 1451.5     & 187.2    & {\bf 46.2}    & 26437.2   & /   \\ \hline
			\multirow{2}{*}{\bf Algorithm \ref{alg6}}  & IT  & 88000    & 236240     & 338280     & 975600     & {\bf 905200}      \\ \cline{2-7}
			& CPU & {\bf 183.4}     & {\bf 72.7}    & 94.6    & {\bf 4393.7}    & {\bf 44789.7}   \\ \hline
		    \multicolumn{2}{|c|}{\bf Speed-up}   &2.3   &4.3   &1.6   &6.9   &/    \\ \hline
	\end{tabular}}\label{T5}}
\end{table}

From Table \ref{T5}, we observe that Algorithm \ref{alg6} is the fastest one while Algorithm \ref{alg5} needs the fewest number of iterations in most cases.
Specifically, Algorithm \ref{alg6} is about two to six times faster than GRAK, and the improvement is prohibitive. For the extremely large-scale data set
$YouTubeFaceArrange\_64\times64$, all the algorithms fail to converge except for Algorithm \ref{alg6}. Thus, we benefit from the strategy of random sampling and there is no need to calculate the whole residual vector at each iteration. All these results show that Algorithm \ref{alg6} is advantageous in solving large-scale least squares problems with dense data. Notice that the PREK method does not work for the dataset $MNIST$, this is because this method may fail to converge if the coefficient matrix is not of full column rank \cite{ref13}.


\section{Concluding Remarks}
The greedy randomized augmented Kaczmarz method is an effective approach proposed recently for large and sparse inconsistent linear systems. The key of this method is to equivalently transform a inconsistent linear system of size $m$-by-$n$ into a consistent augmented linear system of size $(m+n)$-by-$(m+n)$, and then apply the GRK method to the augmented linear system. However, in the GRAK method, it is required to calculate the residual of the augmented linear system and construct two index sets at each iteration. Moreover, it only updates the vector ${\bf z}_k$ rather than the approximate solution ${\bf x}_k$ as $m+1\leqslant t_k\leqslant m+n$, which may slow down the convergence rate. Thus, the computational overhead of this method is large, especially for big-data problems, and it is urgent to enhance the numerical performance of this method.

The main contributions of this work are two-fold. First, we propose an accelerated greedy randomized augmented Kaczmarz method for large-scale and sparse or dense inconsistent linear systems, and apply the strategy of simple random sampling to further reduce the workload. Theoretical analysis shows that, under very weak assumptions, the accelerated greedy randomized augmented Kaczmarz method converges faster than the greedy randomized augmented Kaczmarz method. Second, as far as we know, there are no practical stopping criteria for randomized Kaczmarz-type methods till now. So as to fill-in this gap, we introduce a practical stopping criterion that is applicable to {\it all} the randomized Kaczmarz-type methods. Numerical experiments demonstrate the numerical behavior of our proposed algorithms and the effectiveness of the stopping criterion. However, there is a parameter $L$ involved in the proposed stopping criterion. How to pick its optimal value is interesting and deserves further investigation. We believe it is problem-and-data dependent.

%


\begin{thebibliography}{99}
\bibitem{ref20} {\sc R. Ansorge}, {\em Connections between the Cimmino-method and the Kaczmarz-method for the solution of singular and regular systems of equations}, Computing, 33 (1984), pp. 367--375.
\bibitem{ref15} {\sc V. Borkar, N. Karamchandani, and S. Mirani}, {\em Randomized Kaczmarz for rank aggregation from pairwise comparisons}, IEEE Information Theory Workshop (ITW), Cambridge, 2016, pp. 389--393.
\bibitem{ref3} {\sc Z.-Z. Bai and W.-T. Wu}, {\em On greedy randomized Kaczmarz method for solving large sparse linear systems}, SIAM Journal on Scientific Computing, 40 (2018), pp. A592--A606.
\bibitem{ref6} {\sc Z.-Z. Bai and W.-T. Wu}, {\em On greedy randomized augmented Kaczmarz method for solving large sparse inconsistent linear systems}, SIAM Journal on Scientific Computing, 43 (2021), pp. A3892--A3911.
\bibitem{ref13} {\sc Z.-Z. Bai and W.-T. Wu}, {\em On partially randomized extended Kaczmarz method for solving large sparse overdetermined inconsistent linear systems}, Linear Algebra and its Applications, 578 (2019), pp. 225--250.
\bibitem{ref29} {\sc Z.-Z. Bai and W.-T. Wu}, {\em On relaxed greedy randomized Kaczmarz methods for solving large sparse linear systems}, Applied Mathematics Letters, 83 (2018), pp. 21--26.
\bibitem{ref12} {\sc Z.-Z. Bai, L. Wang, and G.V. Muuratova}, {\em On relaxed greedy randomized augmented Kaczmarz methods for solving large sparse inconsistent linear systems}, East Asian Journal on Applied Mathematics, (27) 2021, pp. 323--332.
\bibitem{ref51} {\sc L. Bai, J. Liang and Y. Zhao}, {\em Self-constrained spectral clustering}, IEEE Transactions on Pattern Analysis and Machine Intelligence, 45 (2023), pp. 5126--5138.
\bibitem{ref38} {\sc J. Chen and Z.-D. Huang}, {\em On a fast deterministic block Kaczmarz method for solving large scale linear system}, Numerical Algorithms, 89 (2022), pp. 1007--1029.
\bibitem{ref11} {\sc S. Chapra and R Canale}, {\em Numerical Methods for Engineers, 7th. ed.}, McGraw-Hill, New York, 2005.
\bibitem{ref45} {\sc Y. Censor}, {\em Row-action methods for huge and sparse systems and their applications}, SIAM Review, 23 (1981), pp. 444--466.
\bibitem{SRDA} {\sc D. Cai, X.-F. He, and J.-W. Han}, {\em SRDA: An Efficient Algorithm for Large-Scale Discriminant Analysis}, IEEE Transactions on Knowledge and Data Engineering, 20 (2008), pp. 1--12.

\bibitem{ref33} {\sc K. Du, W. Si, and X. Sun}, {\em Randomized extended average block Kaczmarz for solving least squares}, SIAM Journal on Scientific Computing, 42 (2020), pp. A3541--A3559.


\bibitem{ref25} {\sc C. Gu and Y. Liu}, {\em Variant of greedy randomized Kaczmarz for ridge regression}, Applied Numerical Mathematics, 143 (2019), pp. 223--246.
\bibitem{ref46} {\sc R. Gordon, R. Bender, and G.T. Herman}, {\em Algebraic reconstruction techniques (ART) for three-dimensional electron microscopy and X-ray photography}, Journal of Theoretical Biology, 29 (1970), pp. 471--481.
\bibitem{ref16} {\sc W. Guo, H. Chen, W. Geng, and L. Lei}, {\em A modified Kaczmarz algorithm for computerized tomographic image reconstruction}, IEEE International Conference on Biomedical Engineering and Informatics, 3 (2009), pp. 1--4.
\bibitem{ref43} {\sc A. Hefny, D. Needell, and A. Ramdas}, {\em Rows versus columns: randomized Kaczmarz or Gauss-Seidel for ridge egression}, SIAM Journal on Scientific Computing, 39 (2017), pp. S528--S542.
\bibitem{ref47} {\sc G. Herman}, {\em Fundamentals of Computerized Tomography: Image Reconstruction from Projections, 2nd. ed.}, Springer, Dordrecht, 2009.
\bibitem{ref14} {\sc P. Hansen and J$\varnothing$rgensen}, {\em AIR tools $\uppercase\expandafter{\romannumeral2}$:algebraic iterative reconstruction methods, improved mplementation}, Numerical Algorithms, 79 (2018), pp. 107--137.

\bibitem{ref37} {\sc X. Jiang, K. Zhang, and J. Yin}, {\em Randomized block Kaczmarz methods with k-means clustering for solving large linear systems}, Journal of Computational and Applied Mathematics, 403 (2022), Article 113828.
\bibitem{ref4} {\sc Y. Jiang, G. Wu, and L. Jiang}, {\em A semi-randomized Kaczmarz method with simple random sampling for large-scale linear systems}, Advances in Computational Mathmatics, 49 (2023), Article 20.
\bibitem{WLaw} {\sc J. Johnson}, {\em Probability and Statistics for Computer Scientists}, The John Wiley \& Sons Press, 1997.    
    
\bibitem{ref1} {\sc S. Kaczmarz}, {\em Approximate solution of systems of linear equations}, International Journal of Control, 35 (1937), pp. 355--357.
\bibitem{ref17} {\sc S. Lee and H. Kim}, {\em Noise properties of reconstructed images in a kilo-voltage on-board imaging system with iterative reconstruction techniques: A phantom study}, Physica Medica, 30 (2014), pp. 365--373.
\bibitem{ref9} {\sc K. Li and G. Wu}, {\em Randomized approximate class-specific kernel spectral regression analysis for large-scale face verification}, Machine Learning, 111 (2022), pp. 2037--2091.
\bibitem{ref35} {\sc R. Li and H. Liu}, {\em  On global randomized block Kaczmarz method for image reconstruction}, Electronic Research Archive, 30 (2022), pp. 1442--1453.
\bibitem{ref23} {\sc J. Liu and S. Wright}, {\em An accelerated randomized Kaczmarz algorithm}, Mathematics of Computation, 85 (2016), pp. 153--178.
\bibitem{ref22} {\sc J. Loera, J. Haddock, and D. Needell}, {\em A sampling Kaczmarz-Motzkin algorithm for linear feasibility}, SIAM Journal on Scientific Computing, 39 (2017), pp. S66--S87.
\bibitem{ref7} {\sc A. Ma, D. Needell, and A. Ramdas}, {\em Convergence properties of the randomized extended Gauss-Seidel and Kaczmarz methods}, SIAM Journal on Matrix Analysis and Applications, 36 (2015), pp. 1590--1604.
\bibitem{ref36} {\sc C.-Q. Miao and W.-T. Wu}, {\em On greedy randomized average block Kaczmarz method for solving large linear systems}, Journal of Computational and Applied Mathematics, 413 (2022), Article 114372.
\bibitem{ref44} {\sc D. Needell}, {\em Randomized Kaczmarz solver for noisy linear systems}, Bit Numerical Mathematics, 50 (2010), pp. 395--403.
\bibitem{ref31} {\sc D. Needell and J. Tropp}, {\em  Paved with good intentions: Analysis of a randomized block Kaczmarz method}, Linear Algebra and Its Applications, 441 (2014), pp. 199--221.
\bibitem{ref34} {\sc D. Needell, R. Zhao, and A. Zouzias}, {\em Randomized block Kaczmarz method with projection for solving least squares}, Linear Algebra and Its Applications, 484 (2015), pp. 322--343.
\bibitem{ref30} {\sc I. Necoara}, {\em Faster randomized block Kaczmarz algorithms}, SIAM Journal on Matrix Analysis and Applications, 40 (2019), pp. 1425--1452.
\bibitem{ref32} {\sc Y. Niu and B. Zheng}, {\em A greedy block Kaczmarz algorithm for solving large-scale linear systems}, Applied Mathematics Letters, 104 (2020), Article 106294.

\bibitem{ref39} {\sc C. Popa}, {\em Least-squares solution of overdetermined inconsistent linear systems using kaczmarz's relaxation}, International Journal of Computer Mathematics, 55 (1995), pp. 79--89.
\bibitem{ref40} {\sc C. Popa}, {\em Extensions of block-projections methods with relaxation parameters to inconsistent and rank-deficient least-squares problems}, Bit Numerical Mathematics, 38 (1998), pp. 151--176.

\bibitem{ref18} {\sc R. Ramlau and M. Rosensteiner}, {\em An efficient solution to the atmospheric turbulence tomography problem using Kaczmarz iteration}, Inverse Problems, 28 (2012), Article 095004.
\bibitem{ref50} {\sc W. Shi and G. Wu}, {\em New algorithms for trace-ratio problem with application to high-dimension and large-sample data dimensionality reduction}, Machine Learning (2021), pp. 1--28.
\bibitem{ref24} {\sc S. Steinerberger}, {\em A weighted randomized Kaczmarz method for solving linear systems}, Mathematics of Computation, 90 (2021), pp. 2815--2826.
\bibitem{ref2} {\sc T. Strohmer and R. Vershynin}, {\em A randomized Kaczmarz algorithm with exponential convergence}, Journal of Fourier Analysis and Applications, 15 (2009), pp. 262--278.

\bibitem{ref19} {\sc G. Thoppe, V. Borkar, and D. Manjunath}, {\em A stochastic Kaczmarz algorithm for network tomography}, Automatica, 50 (2014), pp. 910--914.



\bibitem{ref8} {\sc X. Yang}, {\em A geometric probability randomized Kaczmarz method for large scale linear systems}, Applied Numerical Mathematics, 164 (2021), pp. 139-160.
\bibitem{ref41} {\sc J. Yin, N. Li, and N. Zheng}, {\em Restarted randomized surrounding methods for solving large linear equations}, Applied mathematics letters, 133 (2022), Article 108290


\bibitem{ref27} {\sc J. Zhang}, {\em A new greedy Kaczmarz algorithm for the solution of very large linear systems}, Applied Mathematics Letters, 91 (2019), pp. 207--212.
\bibitem{ref28} {\sc Y. Zhang and H. Li}, {\em Greedy Motzkin-Kaczmarz methods for solving linear systems}, Numerical Linear Algebra with Applications, 29 (2022), e2429.
\bibitem{ref42} {\sc Y. Zhang and H. Li}, {\em A count sketch maximal weighted residual Kaczmarz method for solving highly overdetermined linear systems}, Applied Mathematics and Computation, 410 (2021), Article 126486.
\bibitem{ref5} {\sc A. Zouzias and N.M. Freris}, {\em Randomized extended Kaczmarz for solving least-squares}, SIAM Journal on Matrix Analysis and Applications, 34 (2013), pp. 773--793.


\end{thebibliography}
\end{document}